\documentclass[dvips,aos,preprint]{imsart}

\RequirePackage[OT1]{fontenc} \RequirePackage{amsthm,amsmath,color}
\RequirePackage[numbers]{natbib}
\RequirePackage[colorlinks,citecolor=blue,urlcolor=blue]{hyperref}

\arxiv{}

\startlocaldefs
\numberwithin{equation}{section}
\theoremstyle{plain}

\endlocaldefs

\usepackage[english]{babel}

\usepackage{amssymb,amsfonts,graphics,graphicx,color}

\def\Var{\mathrm{Var}}

\newtheorem{Prop}{Proposition}
\newtheorem{Lemma}{Lemma}
\newtheorem{Hyp}{Assumption}
\newtheorem{Cor}{Corollary}
\newtheorem{Th}{Theorem}

\newcommand{\pa}[1]{\left({#1}\right)}
\newcommand{\cro}[1]{\left[{#1}\right]}

\newcommand{\ac}[1]{\left\{{#1}\right\}}

\newcommand{\calS}{\mathcal{S}}
\newcommand{\norm}[1]{\ensuremath{\vert\!\vert #1 \vert\!\vert}}
\newcommand{\E}{\ensuremath{\mathbb{E}}}

\renewcommand{\P}{\ensuremath{\mathbb{P}}}
\newcommand{\X}{\ensuremath{\mathbb{X}}}
\newcommand{\N}{\ensuremath{\mathbb{N}}}

\newcommand{\R}{\ensuremath{\mathbb{R}}}

\newcommand{\e}{\ensuremath{\varepsilon}}
\newcommand{\la}{\ensuremath{\lambda}}
\newcommand{\La}{\ensuremath{\Lambda}}
\newcommand{\p}{\ensuremath{\varphi}}

\newcommand{\M}{\mathcal{M}}
\newcommand{\1}{{\bf 1}}
\newcommand{\LL}{\mathbb{L}}

\numberwithin{equation}{section}

\begin{document}

\begin{frontmatter}
\title{The two-sample problem for Poisson processes: adaptive tests with a non-asymptotic wild bootstrap approach}
\runtitle{Two-sample problem for Poisson processes}
%\thankstext{T1}{Footnote to the title with the ``thankstext'' command.}

\begin{aug}
\author{\fnms{Magalie} \snm{Fromont}\thanksref{m1}\ead[label=e1]{mfromont@ensai.fr}},
\author{\fnms{B\'eatrice} \snm{Laurent}\thanksref{m2}\ead[label=e2]{Beatrice.Laurent@insa-toulouse.fr}}
\and
\author{\fnms{Patricia} \snm{Reynaud-Bouret}\thanksref{m3}\ead[label=e3]{Patricia.Reynaud-Bouret@unice.fr}}

%\thankstext{t1}{Some comment}
%\thankstext{t2}{First supporter of the project}
%\thankstext{t3}{Second supporter of the project}
\runauthor{M. Fromont et al.}

\affiliation{\thanksmark{m1}CREST (Ensai) - IRMAR (Universit\'e Europ\'eenne de Bretagne), \thanksmark{m2}
IMT, INSA Toulouse, Universit\'e de Toulouse and \thanksmark{m3} Laboratoire Jean-Alexandre Dieudonn\'e, CNRS, Universit{\'e} de Nice
Sophia-Antipolis}

\address{Magalie Fromont\\
ENSAI - Campus de Ker-Lann\\
Rue Blaise Pascal - BP 37203\\
35172 Bruz Cedex (FRANCE)\\
\printead{e1}}

\address{B\'eatrice Laurent\\
D\'epartement GMM -INSA\\
135 av. de Rangueil \\
31077 Toulouse Cedex 4 (FRANCE)\\
\printead{e2}}

\address{ Patricia Reynaud-Bouret\\
Laboratoire Jean-Alexandre Dieudonn\'e\\
CNRS - Universit\'e de Nice Sophia-Antipolis\\
Parc Valrose\\
06108 Nice Cedex 2 (FRANCE)\\
\printead{e3}}
\end{aug}

\begin{abstract} $\ $Considering two independent Poisson processes, we address the
question of testing equality of their respective intensities. We first propose single tests whose test statistics are $U$-statistics based on general kernel functions. The corresponding critical values are constructed from a non-asymptotic wild bootstrap approach, leading to level $\alpha$ tests. Various choices for the kernel functions are possible, including projection, approximation or reproducing kernels. In this last case, we obtain a parametric rate of testing for a weak metric defined in the RKHS associated with the considered reproducing kernel. Then we introduce, in the other cases, an aggregation procedure, which allows us to import ideas coming from
 model selection,
thresholding and/or approximation kernels adaptive estimation.  The resulting multiple tests are proved to be of level
$\alpha$, and to satisfy non-asymptotic oracle type conditions for the classical $\LL_2$-norm.
From these conditions, we deduce that they are
adaptive in the minimax sense over a large variety of classes of
alternatives based on classical and weak Besov bodies in the
univariate case, but also Sobolev and anisotropic Nikol'skii-Besov
balls in the multivariate case.
\end{abstract}



\end{frontmatter}

\section{Introduction}

We consider the two-sample problem  for general Poisson processes.
Let $N^1$ and $N^{-1}$ be two independent Poisson processes observed
on a measurable space $\X$, whose intensities with respect to some
non-atomic positive $\sigma$-finite measure $\mu$ on $\X$ are denoted by $f$
and $g$. Given the observation of $N^1$ and $N^{-1}$, we address the
question of testing the null hypothesis $(H_0)$ "$f=g$" versus the
alternative $(H_1)$ "$f\neq g$".

Many papers deal with the two-sample problem for homogeneous Poisson
processes such as, among others, the historical ones of
\cite{PrzyborowskiWilenski40}, \cite{Cox53}, \cite{Gail74}, or
\cite{ShiueBain82}, whose applications were mainly turned to biology
and medicine, and less frequently to reliability. More recent papers
like \cite{KrishnamoorthyThomson2004}, \cite{NgGuTang2007},
\cite{ChiuWang2009}, and \cite{Chiu2010} give interesting numerical
comparisons of various testing procedures. As for non-homogeneous
Poisson processes, though
a lot of references on the
problem of testing proportionality of the hazard rates of the
processes exist (see \cite{DeshpandeSengupta95} for instance and the
references therein),
 very few papers are devoted to a comparison of the
intensities themselves. Bovett and Saw \cite{BovettSaw80} and
Deshpande et al. \cite{Deshpandeetal99} respectively proposed
conditional and unconditional procedures to test the null hypothesis
"$f/g$ is constant" versus "it is increasing". Deshpande et al.
\cite{Deshpandeetal99} considered their test from a usual asymptotic
point of view, proving that it is consistent against several large
classes of alternatives.

We propose in this paper to construct testing procedures of $(H_0)$
versus $(H_1)$ without any parametric or monotony assumption on $f$ or $g$ and which satisfy specific non-asymptotic performance properties.

 In particular, for every $\alpha$ in $[0,1]$, these tests are of level $\alpha$, that is they have a probability of first kind error at most equal to $\alpha$. For special values of $\alpha$, they are even of size $\alpha$, that is their probability of first kind error is exactly equal to $\alpha$, since they involve very sharp critical values obtained via a non-asymptotic wild bootstrap approach. In the classical two-sample problem for i.i.d. samples, the choice of the critical values in testing procedures
is a well-known crucial question. Indeed, the asymptotic
distributions of many test statistics are not free from the
common unknown density under the null hypothesis. In such cases,
some bootstrap methods are often used to build data-driven
critical values. By bootstrap methods, we mean the original ones introduced by Efron \cite{Efron79} of course, but also more general weighted bootstrap approaches such as the precursor Fisher's \cite{Fisher35} permutation, the $m$ out of
$n$ bootstrap introduced by Bretagnolle \cite{Bretagnolle}, the general exchangeably weighted bootstrap studied in
\cite{PraestgaardWellner} and including the
Bayesian bootstrap of Rubin \cite{Rubin81} for instance, as well as the wild bootstrap detailed in \cite{Mammen92}. Except in the cases where the permutation
approach is used, authors generally prove that the obtained tests are
(only) asymptotically of level $\alpha$ (see among many other papers \cite{Romano88}, \cite{Romano89}, \cite{Praestgaard95}, and more recently
\cite{JanssenPauls} for a complete and very interesting discussion). In this work, we adopt one of these general weighted bootstrap approaches, but
from a non-asymptotic point of view. The critical values of our
tests are constructed from wild bootstrapped $U$-statistics, which are based on Rademacher variables. The use of Rademacher variables is well-known in the bootstrap community since the work of Mammen \cite{Mammen92}, but also particularly in the statistical learning community since the works of Koltchinskii \cite{Koltchinskii2001} and Bartlett et al. \cite{BBL2002}, followed by \cite{Koltchinskii2006}. It was notably proposed for the construction of general confidence bands in a recent paper by Lounici and Nickl \cite{LouniciNickl2011}. The main particularity of our study, as compared with previous ones, is that we prove here that, under $(H_0)$, given the data, the considered wild bootstrapped $U$-statistics exactly have the same distributions as our test statistics. The corresponding tests are consequently of level $\alpha$ for every $\alpha$ in $[0,1]$, and even of size $\alpha$ for particular values of $\alpha$. Note that as in \cite{RW} or in \cite{Hoeffding52}, it is also possible to randomize these tests in order to turn them into size $\alpha$ tests for every $\alpha$. In this sense, our bootstrap method can be
viewed as an adapted version of the permutation bootstrap method in
a Poisson framework. As usual even when permutation methods are considered, the wild bootstrapped critical values of our tests are not computed exactly in practice, but just approximated through a Monte Carlo method. We also address this question from a non-asymptotic point of view, since we also focus on controlling the loss due to the Monte Carlo approximation.

Our test statistics are based on a single kernel function which can be chosen either as a projection kernel, or as an approximation kernel, or as a reproducing kernel.
A non-asymptotic study of the second kind error of our tests is also performed. Given any $\beta$ in $[0,1]$, depending on the chosen kernel, we obtain non-asymptotic conditions which guarantee that the probability of second kind error is at most equal to $\beta$.
This can be done via a sharp control of the wild bootstrapped critical values under the alternative, which results from concentration inequalities for Rademacher chaoses \cite{PenaGine,Latala}.

In order to deduce from
these conditions recognizable asymptotic rates of testing, we assume that
the measure $\mu$ on $\X$ satisfies $d\mu= n d\nu$, where $n$ can be
seen as a growing number whereas the measure $\nu $ is held fixed.
Typically, $n$ may be an integer and the above assumption amounts to
considering the Poisson processes $N^1$ and $N^{-1}$ as $n$ pooled
i.i.d. Poisson processes with respective intensities $f$ and $g$
w.r.t. $\nu$. The reader may also assume for sake of simplicity that $\X$ is a measurable subset of $\R^d$ and that $\nu$ is the
Lebesgue measure, but it is not required: $\nu$
 may be any non-atomic positive $\sigma$-finite measure on any measurable set $\X$. With this normalization, when a reproducing kernel is considered, we obtain a parametric rate of testing for a weak metric defined in the associated RKHS, in the spirit of \cite{Wellner79} or \cite{Gine75} for more classical weak metrics in i.i.d. samples frameworks. Our results complete those of Gretton et al. \cite{Grettonetal}, who introduced reproducing kernels in the two-sample problem for i.i.d. samples. When a projection or an approximation kernel is considered, we obtain the following condition: the probability of second kind error of the test is at most equal to $\beta$ as soon as the $\LL_2$-distance w.r.t. $\nu$ between $f$ and $g$ is larger than a bound, which reproduces a bias-variance decomposition. This bound can be proved to be optimal with an appropriate choice of the vectorial space defining the projection kernel, or of the bandwidth defining the approximation kernel, choice which highly depends on the alternative.

 In order to provide an adaptive test with respect to this choice, we propose to aggregate several of the previous single kernel-based tests, making sure that the resulting multiple test is still of level $\alpha$. We establish oracle type conditions, which guarantee that the probability of second kind error is at most equal to $\beta$. This aggregation approach, inspired by adaptive estimation methods such as model selection, thresholding or approximation kernels methods, was used in many papers devoted to adaptive testing in various classical one-sample frameworks (see \cite{Spokoiny96} or \cite{Spokoiny98} for adaptive tests related to thresholding methods, \cite{Ingster2000} for adaptive tests related to model selection methods, \cite{HorowitzSpokoiny2001} for adaptive tests related to approximation kernels methods, or \cite{BHL} for adaptive tests related to both model selection and thresholding methods for instance). In a Poisson process framework, we proposed in \cite{MBP} an aggregated test of homogeneity also based on both model selection and thresholding approaches. In the two-sample problem for i.i.d. samples, which is closely related to the present problem, Butucea and Tribouley \cite{ButuceaTribouley} propose an adaptive test based on a thresholding approach.

We complete the study by proving that our aggregated tests are also adaptive in a non-asymptotic minimax sense over various classes $\calS_\delta$ of alternatives $(f,g)$ for which $(f-g)$ is smooth with parameter $\delta$. For clarity's sake, let us here recall a few definitions. For any level $\alpha$ test $\Phi_\alpha$, with values in $\{0,1\}$ (rejecting
$(H_0)$ when $\Phi_{\alpha}=1$), one defines its uniform separation rate $\rho(\Phi_\alpha,\mathcal{S}_\delta,\beta)$ over $\mathcal{S}_\delta$ as
\begin{eqnarray}\label{defrate}
\rho(\Phi_\alpha,\calS_\delta,\beta)&=&\inf{}\left\{\rho>0,\sup_{(f,g)\in\calS_\delta,
\|f-g\|> \rho} \P_{f,g}(\Phi_{\alpha}=0)\leq \beta\right\},
\end{eqnarray}
where $\|f-g\|^2=\int (f-g)^2 d\nu$, and $\P_{f,g}$  denotes the joint distribution of $(N^1,N^{-1})$. A level $\alpha$ test $\Phi_\alpha$ is said to be minimax over a particular class $\mathcal{S}_\delta$ if its uniform separation rate achieves its best possible value over $\mathcal{S}_\delta$, which is called the minimax separation rate over $\mathcal{S}_\delta$ (see \cite{Yannick}) up to a multiplicative factor. It is said to be minimax adaptive if its uniform separation rates achieve (up to a possible unavoidable small loss) the minimax separation rates over several classes $\mathcal{S}_\delta$ simultaneously. A great number of papers deal with the computation of the minimax separation rates over various classes of alternatives, or more precisely with the computation of their asymptotic equivalents, that are the minimax rates of testing defined in the key series of papers due to Ingster \cite{Ingster}. The question of the minimax adaptivity has also been widely studied since the work of Spokoiny \cite{Spokoiny96}, who first brought out a context where minimax adaptive testing without a small loss of efficiency is impossible. For the problem of testing the goodness-of-fit of a Poisson process, Ingster and Kutoyants \cite{IngsterKut2007} derived the minimax rate of testing over a Sobolev or a Besov ball. For the problem of testing the homogeneity of a Poisson process, we derived in \cite{MBP} similar minimax results considering classical Besov bodies, and we moreover obtained new minimax adaptivity results considering weak Besov bodies.

In the present two-sample problem for Poisson processes, no previous minimax result is available to our knowledge. As in \cite{MBP}, we here prove that the aggregation of single projection kernel-based tests lead to minimax adaptive tests over some classes of alternatives for which $(f-g)$ belongs to a Besov or a weak Besov body. Such a result can be linked to the minimax results obtained by Butucea and Tribouley \cite{ButuceaTribouley}, noting however that the classes of alternatives they consider impose both $f$ and $g$ to belong to a Besov space, which is more restrictive than only imposing some regularity assumptions on $(f-g)$. Then, when considering the aggregation of single approximation kernel-based tests, we obtain upper bounds for the uniform separation rates over some classes of alternatives based on multivariate Sobolev or anisotropic Nikol'skii-Besov balls. These upper bounds, which are conjectured to be optimal from results of Horowitz and Spokoiny \cite{HorowitzSpokoiny2001} or Ingster and Stepanova \cite{IngsterStepanova2011} in other frameworks, are completely new in our Poisson setting, and even in a general setting for anisotropic Nikol'skii-Besov balls.

The paper is organized as follows. In Section
\ref{single}, we introduce our single kernel-based tests. As explained above, the corresponding critical values are constructed from a wild bootstrap approach, leading to level $\alpha$ single tests. We then give conditions ensuring that
these single tests also have a probability of second kind error at
most equal to $\beta$, and we study the cost due to the
Monte Carlo approximation of the wild bootstrapped critical values.
In Section \ref{multiple}, we construct  level $\alpha$
multiple tests by aggregating several of the single
 tests introduced in Section \ref{single}. Oracle type conditions are obtained, ensuring that these multiple tests have a probability of second kind error at most equal to $\beta$. From these conditions, some of our tests are also proved to be minimax adaptive over various classes of alternatives based on classical and weak Besov bodies in the univariate case, or Sobolev and anistropic Nikol'skii-Besov balls in the multivariate case.
The major proofs are given in Section \ref{preuves}, whereas a simulation study and the other proofs can be found in supplementary materials.

\smallskip

Let us now introduce some notations that will be used all along the
paper. For any measurable function $h$, let when they exist:
$\norm{h}_{\infty}= \sup_{x \in \X} |h(x)|$, and $\norm{h}_1=\int_{\X} |h(x)| d\nu_x$. Recalling that $\norm{h}=(\int_{\X} h(x)^2 d\nu_x)^{1/2}$, we introduce
 the scalar product $\langle.,.\rangle$ associated
with $\norm{.}$. We denote by $dN^1$ and $dN^{-1}$ the point
measures associated with $N^1$ and
 $N^{-1}$ respectively, and to suit for the notation $\P_{f,g}$ of the joint distribution of $(N^1,N^{-1})$, $\E_{f,g}$ stands for the corresponding expectation. We set for any event $\mathcal{A}$ based on $(N^1,N^{-1})$, $\P_{(H_0)}(\mathcal{A})=\sup_{\{(f,g),\ f=g\}} \P_{f,g}(\mathcal{A})$.

Furthermore, we will introduce some constants, that we do not intend to evaluate here, and that are denoted by $C(\alpha,\beta,\ldots)$ meaning that
they may depend on $\alpha$, $\beta$, $\ldots$. Though they are denoted in the same way, they may vary from one line to another.

\smallskip

Finally, let us make the two following assumptions, which together
imply that $f$ and $g$ belong to
 $\LL^{2}(\X,d\nu)$, and which will be satisfied all along the paper, except when specified.

\begin{Hyp}\label{Hypo1}
$\|f\|_{1}<+\infty$ and $\|g\|_{1}<+\infty$.
\end{Hyp}

\begin{Hyp}\label{Hypo2}
$\norm{f}_{\infty}<+\infty$ and $\norm{g}_{\infty}<+\infty$.
\end{Hyp}

\section{Single kernel-based tests with non-asymptotic wild bootstrapped critical values}\label{single}

\subsection{Single kernel-based test statistics}\label{introker}

Since $f$ and $g$ are assumed to satisfy Assumptions \ref{Hypo1} and \ref{Hypo2}, they are also assumed to belong to $\LL^{2}(\X,d\nu)$. Hence, testing $(H_0)$ "$f=g$"  versus $(H_1)$
"$f\neq g$" here amounts to testing that
 "$\|f-g\|=0$" versus "$\|f-g\|>0$".  Considering a well-chosen finite dimensional subspace
$S$ of $\LL^{2}(\X,d\nu)$, if $\Pi_S$ denotes the orthogonal projection onto $S$ for
$\langle.,.\rangle$, any estimator of an increasing function of $\norm{\Pi_S(f-g)}^2$  may thus be a relevant candidate to be a test statistic.
Let $\{\p_\la, \la \in
\La\}$ be an orthonormal basis of $S$ for $\langle.,.\rangle$, and let
$$\hat T= \sum_{\la\in \La} \left(\left(\int_{\X} \p_\la dN^1-\int_{\X} \p_\la dN^{-1}\right)^2 - \int_{\X} \p_\la^2 dN\right),$$
where $N$ is the pooled Poisson process whose point measure is given by $dN=dN^1+dN^{-1}$.
Since $\E\left[\left(\int \p_\la dN^1\right)^2\right]=\left(\int \p_\la(x) f(x) d\mu_x\right)^2+\int \p_\la^2(x)  f(x) d\mu_x$,
 and similarly for $\E\left[\left(\int \p_\la dN^{-1} \right)^2\right]$, recalling that $d\mu=nd\nu$, it is
easy to see that $\hat T$ is an unbiased estimator of
$n^2\norm{\Pi_S(f-g)}^2$, and thus also a possible test statistic, whose large values lead to reject $(H_0)$.

Let $(\e_x^0)_{x\in N}$ be the marks of the points from the pooled process $N$, defined by
$\e_x^0=1$ if the point $x$ of $N$ belongs to $N^1$ and
 $\e_x^0=-1$ if the point $x$ of $N$ belongs to $N^{-1}$. Then $\hat
T$ can also be expressed as
$$\hat T=\sum_{\la\in\La}\sum_{x\neq x'\in N} \p_\la(x)\p_\la(x') \e_{x}^0\e_{x'}^0= \sum_{x\neq x'\in N}\left( \sum_{\la\in\La}\p_\la(x)\p_\la(x') \right)\e_{x}^0\e_{x'}^0.$$
Starting from this remark, we can thus generalize the test statistic $\hat T$ by replacing in its expression the function: $(x,x')\in\X^2\mapsto
\sum_{\la\in\La}\p_\la(x)\p_\la(x')\in \R$ by a general kernel function. So, let $K$ be any symmetric kernel function: $\X\times\X\to \R$
satisfying:
\begin{Hyp}\label{Hypo3} $\int_{\X^2} K^2(x,x')(f+g)(x) (f+g)(x')d\nu_xd\nu_{x'}<+\infty.$
\end{Hyp}
Denoting by $\X^{[2]}$ the set $\{(x,x')\in\X^2,\ x\neq x'\}$, we
introduce the statistic
\begin{equation}\label{defTchapeauK}
\hat T_K=\sum_{x\neq x'\in N} K(x,x')\e_{x}^0\e_{x'}^0=\int_{\X^{[2]}} K(x,x')\e_x^0\e_{x'}^0dN_xdN_{x'}.
\end{equation}

Since for every $x$ in $N$, $\E[\e_x^0|N
]=(f(x)-g(x))/(f(x)+g(x))$ (see Proposition \ref{marques} below for instance),
\begin{eqnarray*}
\E_{f,g}[\hat T_K]&=&\E_{f,g}\left[\E\left[\int_{\X^{[2]}} K(x,x')\e_x^0\e_{x'}^0dN_xdN_{x'}\Big|N\right]\right]\\
&=&\E_{f,g}\left[\int_{\X^{[2]}}K(x,x')\frac{f(x)-g(x)}{f(x)+g(x)}\frac{f(x')-g(x')}{f(x')+g(x')}dN_xdN_{x'}\right]\\
&=& \int_{\X^{2}}K(x,x')(f-g)(x)(f-g)(x')d\mu_xd\mu_{x'}\\
&=& n^2 \int_{\X^{2}}K(x,x')(f-g)(x)(f-g)(x')d\nu_xd\nu_{x'}.
\end{eqnarray*}
In the following, we use the notation:
\begin{equation}\label{diamant}
 K\cro{ p}(x')=\int_{\X}K(x,x')p(x)d\nu_x.
\end{equation}
With this notation, $\hat T_K$ is then an unbiased estimator of
\begin{equation}\label{defEK}
\mathcal{E}_K=n^2\langle K\cro{f-g},f-g\rangle,
\end{equation}
whose existence is ensured thanks to Assumptions \ref{Hypo1} and
\ref{Hypo3}.

\smallskip

We have chosen to consider and study in this paper three possible examples of kernel functions. For each example, we give a simpler expression of $\mathcal{E}_K$, which allows to justify the choice of $\hat T_K$ as test statistic.

\smallskip

\emph{[Projection kernel case]} Our first choice for $K$ is a symmetric kernel function based on
an orthonormal family $\{\p_\la, \la \in \La\}$ for
$\langle.,.\rangle$:
$$K(x,x')=\sum_{\la\in\La}\p_\la(x)\p_\la(x').$$
When  the cardinality of $\La$ is finite, $\hat{T}_K$ corresponds to
the above natural test statistic $\hat{T}$.  When the cardinality
of $\La$ is infinite, we assume that
$\sup_{x,x'\in \X} \sum_{\la \in \Lambda} |\p_\la(x)\p_\la(x')| <+\infty,$
which ensures that $K(x,x')$ is defined for all $x,x'$ in $\X$ and
that Assumption~\ref{Hypo3} holds. Typically, if $\X=\R^d$ and if
the functions $\{\p_\la, \la \in \La\}$ correspond to indicator
functions with disjoint
supports, this condition will be satisfied.

We check in these cases that for every $s$ in $\LL^2(\X,d\nu)$, $K\cro{s}=\Pi_S(s),$
where $S$ is the subspace of $\LL^2(\X,d\nu)$ generated by  $\{\p_\la, \la \in \La\}$, and $\Pi_S$ denotes as above the orthogonal projection onto
$S$ for $\langle.,.\rangle$. This justifies that such a kernel function $K$ is called a projection kernel and that
$$\mathcal{E}_K=n^2 \norm{\Pi_S(f-g)}^2.$$

\smallskip

\emph{ [Approximation kernel case]} When $\X=\R^d$ and $\nu$ is the Lebesgue measure, our second
choice for $K$ is a kernel function based on an approximation kernel
$k$ in $\LL^2(\R^d)$, and such that $k(-x)=k(x)$: for
$x=(x_1,\ldots,x_d)$, $x'=(x_1',\ldots,x_d')$ in
 $\X$,
$$K(x,x')= \frac{1}{\prod_{i=1}^d h_i}k\pa{\frac{x_1-x_1'}{h_1},\ldots,\frac{x_d-x_d'}{h_d}},$$
where $h=(h_{1},\ldots,h_{d})$ is a vector of $d$ positive
bandwidths. Note that the assumption that $k\in\LL^2(\R^d)$
together with Assumption \ref{Hypo2} ensure that Assumption
\ref{Hypo3} holds. Then, in this case,
$$\mathcal{E}_K= n^2\langle
k_h*(f-g),f-g\rangle,$$ where
$k_h(u_1,\ldots,u_d)=\frac{1}{\prod_{i=1}^d
h_i}k\pa{\frac{u_1}{h_1},\ldots,\frac{u_d}{h_d}}$ and $*$ is the
usual convolution operator with respect to the measure  $\nu$.

\smallskip

\emph{ [Reproducing kernel case]} Our third choice for $K$ is a general reproducing kernel
(see \cite{ScholSmola} for instance) such that
$$K(x,x')=\langle \theta(x),\theta(x')\rangle_{\mathcal{H}_K},$$
where $\theta$ and $\mathcal{H}_K$ are a representation function and
a RKHS associated with $K$. Here, $\langle .,.
\rangle_{\mathcal{H}_K}$ denotes the scalar product of
$\mathcal{H}_K$. We also choose $K$ such that it satisfies
Assumption \ref{Hypo3}.

This choice leads to a test statistic
close to the one of Gretton et al. \cite{Grettonetal} for the classical
two-sample problem for i.i.d. samples of equal sizes.
We will however see that the corresponding critical value is not constructed here in
the same way as in \cite{Grettonetal}. While Gretton et al. derive their critical value from either
concentration inequalities, or asymptotic arguments, or an asymptotic
Efron's bootstrap approach, we construct our critical value from
a non-asymptotic wild bootstrap approach.

In this case, it is easy to see that
$$\mathcal{E}_K=n^2\left\|m_f-m_g\right\|_{\mathcal{H}_K}^2,$$
where $m_f=\int_{\X}  K(.,x)f(x)d\nu_x$ and $m_g=\int_{\X}  K(.,x)g(x)d\nu_x$.
Note that in a "density" context where $\int_\X f(x) d\nu_x =\int_\X g(x) d\nu_x=1$, $\mathcal{E}_K$ is $n^2$ times the so-called squared Maximum Mean Discrepancy on the unit ball
in the RKHS $\mathcal{H}_K$ (see \cite{Grettonetal}) between the distributions $fd\nu$ and $gd\nu$, and that the functions $m_f$ and $m_g$ are known (see \cite{SriFukumizuLanckriet} for instance) as the mean embeddings in $\mathcal{H}_K$ of the distributions $fd\nu$ and $gd\nu$ respectively. Moreover, in this context, assuming that the kernel $K$ is characteristic (see also \cite{SriFukumizuLanckriet}), the map which assigns its mean embedding in $\mathcal{H}_K$ to any probability distribution is injective by definition, so
$\mathcal{E}_K=0$ if and only if $f=g$.

We want to mention here that the introduction of reproducing kernels is particularly pertinent if the space $\X$
is unusual or pretty large
 with respect to the (mean) number of observations and/or if the measure $\nu$ is not well specified or not easy to deal with. In such
situations, the use of reproducing kernels may be the only possible way
to compute a meaningful test (see \cite{Grettonetal} where such
kernels are used for microarrays data and graphs).

\smallskip

Thus, for each of the three above choices for $K$, considering a test which rejects $(H_0)$ when $\hat T_K$ is "large enough" seems to be reasonable. It remains to explain what we mean by "large enough", that is to define the critical values used in our tests.

 \subsection{Critical values based on a non-asymptotic wild bootstrap approach}\label{bstcriticalvalue}

 The critical values we use here are based on a non-asymptotic wild bootstrap approach, that we present and justify in this section.
 To do this, we start from the remark that under $(H_0)$, the test statistic $\hat T_K$ is a degenerate $U$-statistic of order $2$, for which adequate bootstrap methods were developed in particular in \cite{Bretagnolle} and \cite{ArconesGine92}. Bretagnolle \cite{Bretagnolle} first noticed that a naive application of Efron's original bootstrap fails for degenerate $U$-statistics, since it leads the bootstrapped statistic to lose the degeneracy property. He therefore introduced the more appropriate $m$ of $n$ bootstrap, while Arcones and Gin\'e~\cite{ArconesGine92} preferred to keep on using Efron's original bootstrap, but by forcing the bootstrapped statistic to satisfy the degeneracy property through a centering trick. The results of Arcones and Gin\'e were then generalized to other kinds of bootstrap methods, and in particular Bayesian and wild bootstrapped $U$-statistics were introduced in \cite{HuskovaJanssen}, \cite{Janssen94} and \cite{DehlingMikosch}.

Following \cite{DehlingMikosch}, we introduce a sequence $(\e_i)_{i\in\N}$ of
i.i.d. Rademacher variables independent of $N$. Denoting by $N_n$ the size of the pooled process $N$, and by $\{X_1,\ldots, X_{N_n}\}$
the points of $N$, a wild bootstrapped version of $\hat T_K$ may be expressed as
$\sum_{i\not = i' \in \{1,\ldots, N_n\}}K(X_i,X_{i'})\e_{X_i}^0\e_{X_{i'}}^0 \e_i\e_{i'}.$
We consider in fact the simpler version
\begin{equation}\label{defTchapeaubootK}
\hat{T}_K^{\e}=\sum_{i\neq i'\in\{1,\ldots,N_n\}} K(X_i,X_{i'}) \e_i\e_{i'},
\end{equation}
that can be proved to have, under $(H_0)$, conditionally on $N$, the same distribution as the above wild bootstrapped version of
$\hat T_K$.
We now choose the quantile of the conditional distribution of $\hat{T}_K^{\e}$ given $N$ as critical value for our test. 

More precisely, for $\alpha$ in $(0,1)$, if $q^{(N)}_{K,1-\alpha}$ denotes the
$(1-\alpha)$ quantile of the distribution of $\hat{T}_K^{\e}$ conditionally on $N$, we consider the test that rejects $(H_0)$ when $\hat
T_K>q^{(N)}_{K,1-\alpha}$.
The corresponding test function  is defined by
\begin{equation}\label{fonctiontest}
\Phi_{K,\alpha}=\1_{\hat T_K>q^{(N)}_{K,1-\alpha}}.
\end{equation}
Note that in practice, the true conditional quantile $q^{(N)}_{K,1-\alpha}$ is not exactly computed, but in fact just approximated by a classical Monte Carlo method.

Of course, such bootstrap tests are not completely new in the statistical scene. However, the main particularities of our work is that we justify our test from a non-asymptotic point of view. We actually prove that under
$(H_0)$, conditionally on $N$,
$\hat{T}_K$  and $\hat{T}_K^\e$ exactly have the same distribution.
As a
consequence the test defined by $\Phi_{K,\alpha}$ is of level $\alpha$, that is it has a probability of first kind error at most equal to $\alpha$. We will briefly see in the next section that it may even be randomized to be of size $\alpha$, that is to have a probability of first kind error exactly equal to $\alpha$.

In the same
way, instead of focusing as many previous authors on the
consistence against some alternatives, we give precise
conditions on the alternatives which guarantee that $\Phi_{K,\alpha}$  has a probability of second kind error controlled by a prescribed value $\beta$ in $(0,1)$.
These results are detailed in the next section.

Furthermore, we do not forget that studying our tests from a
non-asymptotic point of view poses the additional question
of the exact loss in probabilities of first and second kind errors due to the Monte Carlo approximation of $q^{(N)}_{K,1-\alpha}$. We also address this question in Section
\ref{MonteCarlo}.

Such a non-asymptotic approach is actually conceivable thanks to the following proposition,
which can be deduced from a general result of \cite{daleyVerejones},
but whose quite easy and complete  proof is given in Section
\ref{preuves} for sake of understanding.

\begin{Prop}\label{marques} Let $N^1$ and $N^{-1}$ be two independent Poisson processes on a metric space
$\X$ with intensities $f$ and $g$ with respect
 to some measure $\mu$ on $\X$ and such that Assumption \ref{Hypo1} is satisfied.  Then the pooled process $N$ whose
 point measure is given by $dN=dN^1+dN^{-1}$ is a Poisson process on $\X$ with intensity $f+g$ with respect to $\mu$.
Moreover, let $\left(\e^0_x\right)_{x\in N}$ be defined by  $\e_x^0=1$ if $x$
belongs to $N^1$ and
 $\e_x^0=-1$ if $x$ belongs to $N^{-1}$. Then, conditionally on $N$, the variables $\left(\e^0_x\right)_{x\in N}$
are i.i.d. and for every $x$ in $N$,
\begin{equation}\label{loieps0}
\P\left(\e^0_x=1
|N\right)=\frac{f\left(x\right)}{f\left(x\right)+g\left(x\right)},\
\P\left(\e^0_x=-1
|N\right)=\frac{g\left(x\right)}{f\left(x\right)+g\left(x\right)},
\end{equation}
with the convention that $0/0=1/2$.
\end{Prop}

\subsection{Probabilities of first and second kind errors}\label{secsingleerror}

We here study the probabilities of first and second kind errors of the test $\Phi_{K,\alpha}$ defined by (\ref{fonctiontest}).

\smallskip

From Proposition \ref{marques}, we deduce that under $(H_0)$, $\hat T_K$  and $\hat{T}_K^{\e}$ exactly have the same distribution
 conditionally on $N$. As a result, given $\alpha$ in $(0,1)$, under~$(H_0)$,
\begin{equation}\label{inequationnivcond}
\P\left(\hat T_K>q^{(N)}_{K,1-\alpha}\Big| N\right)\leq \alpha.
\end{equation}
By taking the expectation over $N$, we obtain that
$$\P_{(H_0)}(\Phi_{K,\alpha}=1)\leq \alpha.$$
 In fact, the inequality (\ref{inequationnivcond}) can be turned in an equality only for some particular values of $\alpha$, due to the discreteness of the conditional distribution of $\hat{T}_K$ given $N$.
To go a little further, from Proposition \ref{marques}, we deduce that the randomization hypothesis as defined by Romano and Wolf \cite{RW} and introduced by Hoeffding \cite{Hoeffding52} is satisfied. From the construction of Hoeffding \cite{Hoeffding52}, one can therefore randomize $\Phi_{K,\alpha}$ to obtain a
test $\Psi_{K,\alpha}$ such that $\Psi_{K,\alpha}\geq \Phi_{K,\alpha}$ a.s. and such that under $(H_0)$, $\E(\Psi_{K,\alpha}|N)=\alpha$ for every $\alpha$. Thus, by using the classical tool of randomization, one can circumvent the trouble due to the atoms of the discrete conditional distribution of $\hat{T}_K$ given $N$, and obtain a test with a probability of first kind error exactly equal to $\alpha$ for every $\alpha$. Note that the randomized
test $\Psi_{K,\alpha}$ necessarily has a probability of second kind error smaller than $\Phi_{K,\alpha}$'s one, since $\Psi_{K,\alpha}\geq \Phi_{K,\alpha}$ a.s.

However, in practice, since the conditional quantile $q^{(N)}_{K,1-\alpha}$ is approximated by a Monte Carlo method as we have explained above, we do not have access to the true randomized version of $\Phi_{K,\alpha}$.
This explains why we have decided to focus in the following on the non-randomized test $\Phi_{K,\alpha}$.

\smallskip

Given $\beta$ in $(0,1)$, we now aim at bringing out a non-asymptotic condition on the alternative $(f,g)$ which will guarantee that $\P_{f,g}(\Phi_{K,\alpha}=0)\leq \beta.$
Denoting by $q_{K,1-\beta/2}^\alpha$ the $(1-\beta/2)$ quantile of the conditional quantile
$q^{(N)}_{K,1-\alpha}$,
\begin{equation*}
\P_{f,g}(\Phi_{K,\alpha}=0)\leq \P_{f,g}(\hat{T}_K\leq
q_{K,1-\beta/2}^\alpha) +\beta/2.
\end{equation*}
Thus, a condition which guarantees that $\P_{f,g}(\hat{T}_K\leq q_{K,1-\beta/2}^\alpha)\leq \beta/2$  will
be enough to ensure that $\P_{f,g}(\Phi_{K,\alpha}=0)\leq \beta.$
The following proposition gives such a condition.

\begin{Prop}\label{propsinglerror}
Let $\alpha,\beta$ be fixed levels in $(0,1)$, and let us recall that for any symmetric kernel function $K$ satisfying
Assumption \ref{Hypo3}, $\E_{f,g}[\hat T_K]=\mathcal{E}_K$, with $\mathcal{E}_K$ given in (\ref{defEK}). If
\begin{equation}\label{condek}
\mathcal{E}_K>2n\sqrt{\frac{2nA_K+B_K}{\beta}}+q_{K,1-\beta/2}^\alpha,
\end{equation}
with $A_K=\int_\X\left(K\cro{f-g}(x)\right)^2(f+g)(x)d\nu_x,$ and
$B_K= \int_{\X^2} K^2(x,x')(f+g)(x)(f+g)(x')d\nu_xd\nu_{x'},$
then $\P_{f,g}(\hat{T}_K\leq q_{K,1-\beta/2}^\alpha)\leq \beta/2, $  so that
 $$\P_{f,g}(\Phi_{K,\alpha}=0)\leq \beta.$$
Moreover, there exists some constant $\kappa>0$ such that, for every $K$,
\begin{equation}\label{controlq}
 q_{K,1-\beta/2}^\alpha \leq \kappa \ln (2/\alpha) n \sqrt{\frac{2B_K}{\beta}}.
 \end{equation}
\end{Prop}
To prove the first part of this result, we simply use Markov's inequality since obtaining precise constants and dependency in $\beta$ is not crucial here (see Section \ref{preuves}). The control of $q_{K,1-\beta/2}^\alpha$ derives from a property of Rademacher chaoses combined with an exponential inequality (see \cite{PenaGine} and \cite{Latala}).

\smallskip

The following theorem allows to better understand Proposition \ref{propsinglerror}, and to deduce from it more recognizable properties in terms of uniform separation rates.

\begin{Th}\label{singlerror}
Let $\alpha,\beta$ be fixed levels in $(0,1)$. Let $K$ be a symmetric kernel function satisfying Assumption \ref{Hypo3}, and $\Phi_{K,\alpha}$
be the test defined by~(\ref{fonctiontest}).
Let $C_K$ be an upper bound for $\int_{\X^2} K^2(x,x')(f+g)(x)(f+g)(x')d\nu_xd\nu_{x'}$. Then, we have $\P_{f,g}(\Phi_{K,\alpha}=0)\leq \beta$, as soon as
\begin{multline}\label{RHS}
\norm{f-g}^2\geq
\inf_{r>0}\Big[\left\|(f-g)-r^{-1}K\cro{f-g}\right\|^2\\
+\frac{4+2\sqrt{2}\kappa
\ln(2/\alpha)}{nr\sqrt{\beta}}\sqrt{C_K}\Big]+\frac{8\norm{f+g}_\infty}{\beta
n}.
\end{multline}
For instance, $C_K$ can be taken as follows.
\begin{itemize}
\item $C_K=\|f+g\|_\infty^2 D$ when $K$ is chosen as in the [Projection kernel case], considering an orthonormal basis $\{\p_\la, \la \in \La\}$ of a $D$-dimensional subspace $S$ of $\LL^{2}(\X,d\nu)$,
\item $C_K=\|f+g\|_\infty\|f+g\|_1 D$ when $K$ is chosen as in the [Projection kernel case], considering an orthonormal basis $\{\p_\la, \la \in \La\}$ of a
possibly infinite dimensional subspace $S$ of $\LL^{2}(\X,d\nu)$, which satisfies:
\begin{eqnarray}
&&\sup_{x,x'\in \X} \sum_{\la\in\La}
|\p_\la(x)\p_\la(x')|=D <+\infty, \label{AA1}\\
&&\int_{\X^2} \left( \sum_{\la\in\La}|\p_\la(x)\p_\la(x') |\right)^2
(f+g)(x')  d\nu_xd\nu_{x'} <+\infty,\label{AA2}
\end{eqnarray}
\item  $C_K=\norm{f+g}_\infty\norm{f+g}_{1}\norm{k}^2/
\prod_{i=1}^d h_i$ when $K$ is chosen as in the [Approximation kernel case].
\end{itemize}

\end{Th}

\emph{Comments.}

1. When $K$ is chosen as in the \emph{[Projection kernel case]}, then $K\cro{f-g}=\Pi_S(f-g)$. Hence by taking $r=1$ in (\ref{RHS}), the right hand side of
the inequality reproduces a bias-variance decomposition close to the
bias-variance decomposition for projection estimators, with a
variance term of order $\sqrt{D}/n$ instead of $D/n$. This is quite
usual for this kind of test (see \cite{Yannick} for instance), and
we know that this leads to sharp upper bounds for the uniform
separation rates over particular classes of alternatives.

2. When $K$ is chosen as in the \emph{[Approximation kernel case]} with $k$ in $\LL^1(\R^d)$, $\int_{\R^d} k(x) d\nu_x =1$, and $h_1=\ldots=h_d$,
then $K\cro{f-g}=k_h*(f-g)$, and $\norm{(f-g)-K\cro{f-g}}$ is a bias term. Hence by taking $r=1$ in the inequality (\ref{RHS}), we still reproduce a bias-variance decomposition, but with a variance term of order $h_1^{-d/2}/n$, which coincides with the above variance term in the \emph{[Projection kernel case]} through the equivalence $h_1^{-d} \sim D$. This equivalence is usual in the approximation estimation theory (see
\cite{Sacha} for instance for more details).

3. When $K$ is chosen as in the \emph{[Reproducing kernel case]}, if $K$ is proportional to a kernel from the two above cases, then one can appropriately choose the constant $r$ such that
 $\norm{(f-g)-r^{-1}K\cro{f-g}}$ is still a bias term. We thus recover for such kernel functions, such as the Gaussian and Laplacian kernels, which are commonly used in statistical learning theory, the same bias-variance decomposition as above. However, in some cases, one can not find any normalization constant $r$ for which
$\norm{(f-g)-r^{-1}K\cro{f-g}}$ can be viewed as a
bias term, and the result can not be interpreted from a statistical point of view. In these cases in particular, the $\LL^2$-norm which is considered in Theorem \ref{singlerror} is not the appropriate one to obtain relevant uniform separation rates, since it does not necessarily have any link with the norm of the RKHS $\mathcal{H}_K$. We give in the following theorem a more adequate result for the specific \emph{[Reproducing kernel case]}.

\begin{Th}\label{vitesseparametrique}
Let $\alpha,\beta$ be fixed levels in $(0,1)$, and $\kappa>0$ be the
constant of Proposition \ref{propsinglerror}. Let $\X=\R^d$ and $K$ be a kernel function on $\X\times\X$ chosen as in the [Reproducing kernel case].
Let $\Phi_{K,\alpha}$
be the test function defined by (\ref{fonctiontest}).
We assume furthermore that $\int_\X f(x)d\nu_x=\int_\X g(x)d\nu_x=1$, that $K$ is a bounded measurable characteristic kernel, and that $K(x,x)$ is constant equal to $\kappa_0$. Let $m_f$ and $m_g$ be the mean embeddings of the
distributions $fd\nu$ and $gd\nu$ respectively in $\mathcal{H}_K$. We have $\P_{f,g}(\Phi_{K,\alpha}=0)\leq \beta$ if
$$\|m_f-m_g\|_{\mathcal{H}_K}^2 \geq \frac{4\kappa_0}{n}\pa{ \frac{4 }{\beta} +
\frac{2+\kappa\sqrt{2} \ln(2/\alpha)}{\sqrt{\beta} }}.$$
\end{Th}

\emph{Comments.}

1. The assumption that $K(x,x)$ is constant is usual, since it is satisfied by any normalized or translation-invariant kernel (see \cite{ScholSmola} p 46-47, 57, or
 \cite{SriFukumizuLanckriet} for instance). Moreover, as specified in \cite{SriFukumizuLanckriet} for instance, bounded continuous characteristic and translation-invariant reproducing kernels exist, at least in $\R^d$, where Bochner's theorem enables to characterize them.

2. The result that we have here is in fact comparable to the one obtained by Wellner \cite{Wellner79} for two-sample tests in an i.i.d. samples
 framework. While Wellner's test is based on the estimation of a weak distance between $fd\nu$ and $gd\nu$, associated with the Sobolev norm with
negative index, our test statistic is an unbiased estimator of $\mathcal{E}_K=n^2\|m_f-m_g\|_{\mathcal{H}_K}^2$, where
 $\|m_f-m_g\|_{\mathcal{H}_K}=\sup_{\|r\|_{\mathcal{H}_K}\leq 1} \int_\X(f-g)(x)r(x)d\nu_x$ defines a weak distance between the distributions
$fd\nu$ and $gd\nu$. As in \cite{Wellner79} (or \cite{Gine75} beforehand for the problem of testing uniformity), we obtain a uniform separation
rate for this weak distance of the same order as the usual parametric separation rate, that is of order $n^{-1/2}$.

\subsection{Performance of the Monte Carlo approximation}\label{MonteCarlo}
\subsubsection{Probability of first kind error}

In practice, a Monte Carlo method is used to approximate the conditional  quantiles
 $q^{(N)}_{K,1-\alpha}$. It is therefore necessary to address the following question: what can we say about
the probabilities of first and second kind errors of the test built
with these Monte Carlo approximations?  Recall that we consider the
test $\Phi_{K,\alpha}$ rejecting $(H_0)$ when $\hat{T}_K>
q^{(N)}_{K,1-\alpha}$, where $\hat{T}_K$ is defined by
(\ref{defTchapeauK}), and $q^{(N)}_{K,1-\alpha}$  is the
$(1-\alpha)$ quantile of $\hat{T}_K^\e$ defined by
(\ref{defTchapeaubootK}) conditionally on $N$. The conditional quantile
$q^{(N)}_{K,1-\alpha}$ is estimated by $\hat{q}^{(N)}_{K,1-\alpha}$
via the Monte Carlo method as follows. Conditionally on $N$, we
consider  a set of $B$ independent sequences $\{\e^b, 1\leq b \leq
B\}$, where $ \e^b= (\e^b_x)_{x\in N}$ is a sequence of i.i.d.
Rademacher random variables. We define, for $1\leq b \leq B$,
$\hat{T}_K^{\e^b}= \sum_{x\neq x'\in N } K(x,x')\e_{x}^b\e_{x'}^b.$
Under $(H_0)$, conditionally on $N$, the variables $
\hat{T}_K^{\e^b}$ have the same distribution function as
$\hat{T}_K$, which  is denoted by $F_K$. We denote by $F_{K,B}$ the
empirical distribution function (conditionally on $N$)
 of the sample $( \hat{T}_K^{\e^b},  1\leq b\leq B)$:
$$\forall x \in \R,\  F_{K,B}(x)= \frac{1}{B}\sum_{b=1}^B \1_{ \hat{T}_K^{\e^b} \leq x}.$$
Then, $ \hat{q}^{(N)}_{K,1-\alpha}$  is defined by $ \hat{q}^{(N)}_{K,1-\alpha} = \inf \ac{t \in \R, F_{K,B}(t) \geq 1-\alpha }.$
We finally consider the test given by
\begin{equation}\label{testMC}
 \hat{\Phi}_{K,\alpha}= \1_{\hat{T}_K >\hat{q}^{(N)}_{K,1-\alpha}}.
\end{equation}

\begin{Prop} \label{MC1}
Let $\alpha$ be some fixed level in $(0,1)$, and $ \hat{\Phi}_{K,\alpha}$ be the test defined by (\ref{testMC}).
 Under $(H_0)$,
$$ \P\pa{ \hat{\Phi}_{K,\alpha} =1 \Big|N} \leq \frac{\lfloor B \alpha \rfloor+1}{B+1}.$$
\end{Prop}
\emph{Comment.} For example, if $B=200$ and $\alpha=0.05$, $ \hat{\Phi}_{K,\alpha}$ is of level $5.5\%$.

\subsubsection{Probability of second kind error}

\begin{Prop}\label{puisMC}
Let $\alpha$ and $\beta$ be fixed levels in $(0,1)$ such that $\alpha_B=\alpha-\sqrt{\ln{B}/(2B)}>0$ and $\beta_B=\beta-2/B>0$. Let $ \hat{\Phi}_{K,\alpha}$ be the test given in (\ref{testMC}). Let $\mathcal{E}_K$, $A_K$, $B_K$ and $\kappa$ as in Proposition  \ref{propsinglerror}, and let $q^{\alpha_B}_{K,1-\beta_B/2}$ be the $(1-\beta_B/2)$ quantile of $q_{K,1-\alpha_B}^{(N)}$. If
\begin{equation}
\label{condekMC1}
\mathcal{E}_K > 2n\sqrt{\frac{2nA_K+B_K}{\beta}} + q^{\alpha_B}_{K,1-\beta_B/2},
\end{equation}
then $\P_{f,g}(\hat{\Phi}_{K,\alpha} =0)\leq \beta.$
Moreover,
\begin{equation}
\label{condekMC2}
q^{\alpha_B}_{K,1-\beta_B/2}\leq  \kappa \ln(2/\alpha_B) n\sqrt{\frac{2B_K}{\beta_B}}.
\end{equation}
\end{Prop}
\emph{Comments.} When comparing (\ref{condekMC1}) and (\ref{condekMC2}) with  (\ref{condek})
and (\ref{controlq}) in Proposition~\ref{propsinglerror}, we
notice that they asymptotically coincide when $B\to+\infty$.
Moreover, if $\alpha=\beta=0.05$ and $B\geq 6000$, the
multiplicative factor of $\kappa n\sqrt{B_K}$ is multiplied by a
factor of order $1.2$ in  (\ref{condekMC2}) compared with
(\ref{controlq}). If even $B=200000$, this factor passes
 from $23.4$ in (\ref{controlq}) to $24.1$ in (\ref{condekMC2}).

\section{Multiple testing  procedures}\label{multiple}

In the above section, we  consider testing procedures
based on a single kernel function $K$. Using such single tests
however leads to the natural question of the choice of the kernel,
and/or its parameters: the orthonormal family when $K$ is a projection kernel, the vector of bandwidths $h$ when $K$ is
based on an approximation kernel, the parameters of $K$ when it is a
reproducing kernel. Authors often choose particular parameters regarding
the performance properties that they target for their tests, or use
a data-driven method to choose these parameters which is not always
justified. For instance, in \cite{Grettonetal}, the parameter of the
kernel is chosen from a heuristic method.

 In order to avoid choosing particular kernels or parameters, we propose in this section to consider some collections of kernel
functions instead of a single one, and to define multiple testing
procedures by aggregating the corresponding single tests. We propose
an adapted choice for the critical value. Then, we prove that these
multiple tests satisfy strong statistical properties, such as oracle type properties and minimax adaptivity properties over many classes
of alternatives.

\subsection{Description of the multiple testing procedures}

Let us introduce a finite collection  $\{K_m,m\in \M\}$ of symmetric
kernel functions: $\X\times\X\to \R$  satisfying Assumption
\ref{Hypo3}. For every $m$ in $\M$,  let $\hat{T}_{K_m}$ and  $
\hat{T}_{K_m}^\e $ be defined by (\ref{defTchapeauK})
and (\ref{defTchapeaubootK}) respectively, with $K=K_m$, and let $\{w_m,m\in \M\}$
be a collection of positive numbers such that $\sum_{m\in \M}
e^{-w_m}\leq 1$. For $u$ in $(0,1)$, we denote by $q^{(N)}_{m,1-u}$ the
$(1-u)$ quantile of $\hat{T}_{K_m}^\e$ conditionally on the pooled
process $N$. Given $\alpha$ in $(0,1)$, we consider the test which
rejects $(H_0)$ when there exists at least one $m$ in $\M$ such that
\begin{equation*}
\hat{T}_{K_m}> q^{(N)}_{m,1-u_{\alpha}^{(N)}e^{-w_m}},
\end{equation*}
where $ u_{\alpha}^{(N)}$ is defined by
\begin{equation}\label{ualpha}
u_{\alpha}^{(N)} =\sup \ac{ u > 0, \P\pa{ \sup_{m \in \M}\left(
\hat{T}_{K_m}^\e- q^{(N)}_{m,1-ue^{-w_m}}\right)>0 ~\Bigg| ~N } \leq
\alpha}.
\end{equation}
 Let $\Phi_\alpha$ be the corresponding test
function defined by
\begin{equation}
\label{testmult}
\Phi_\alpha=\1_{\sup_{m\in\M} \left(\hat{T}_{K_m} -q^{(N)}_{m,1-u_{\alpha}^{(N)}e^{-w_m}}\right)>0}.
\end{equation}
Note that given the pooled process $N$,
$u_{\alpha}^{(N)}$ and the quantile $q^{(N)}_{m,1-u_{\alpha}^{(N)}
e^{-w_m}}$ can be estimated by a  Monte Carlo method.

It is quite straightforward to see that this test is of level
$\alpha$ and that one can guarantee a probability of second kind
error at most equal to $\beta$ in $(0,1)$ if one can guarantee it for
one of the single tests rejecting $(H_0)$ when $\hat{T}_{K_m}
>q^{(N)}_{m,1-u_{\alpha}^{(N)}e^{-w_m}}$. We can thus combine the
 results of Theorem~\ref{singlerror}.

\subsection{Oracle type conditions for the probability of second kind error}

\subsubsection{Multiple testing procedures based on projection kernels}\label{sectiontestmulti}

\begin{Th}\label{testmulti}
Let $\alpha,\beta$ be fixed levels in $(0,1)$.
 Let $\{S_m, m\in \M\}$ be a
finite collection of linear subspaces of $\mathbb{L}^2(\X,d\nu)$ and
for all $m$ in $\M$, let $\{\p_\la, \la \in \La_m\}$ be an
orthonormal basis of $S_m$ for $\langle.,.\rangle$.
 We assume either that $S_m$ has finite dimension $D_m$  or that the conditions (\ref{AA1}) and (\ref{AA2}) hold  with $\Lambda=\Lambda_m$ and $D=D_m$.
 We set, for all $m$ in $\M$, $K_m(x,x')= \sum_{\la\in \La_m} \p_\la(x)\p_\la(x') $. Let
$\Phi_{\alpha}$ be the test defined by (\ref{testmult})
with the collection of kernels $\{K_m,m\in\M\}$ and a collection
$\{w_m,m\in \M\}$ of positive numbers such that $\sum_{m\in \M}
e^{-w_m}\leq 1$.

Then $\Phi_{\alpha}$ is a level $\alpha$ test. Moreover,
$\P_{f,g}\pa{\Phi_{\alpha}=0 } \leq \beta$ if
\begin{multline}\label{oraclemulti}
\norm{f-g}^2\geq \inf_{m\in \M}
\Bigg\{\norm{(f-g)-\Pi_{S_m}(f-g)}^2\\
+ \frac{4 + 2 \sqrt{2}
\kappa(\ln(2/\alpha)+w_m)}{n\sqrt{\beta}}M(f,g) \sqrt{D_m}\Bigg\}
+\frac{8\norm{f+g}_\infty}{\beta n},
\end{multline}
where $\kappa>0$ and $M(f,g)= \max \pa{\norm{f+g}_\infty, \sqrt{\norm{f+g}_\infty
\norm{f+g}_1}}$.
\end{Th}

\emph{Comments.} Comparing this result with the one obtained in Theorem \ref{singlerror} for the single test based on a projection kernel, one can see that considering the multiple testing procedure allows
to obtain the infimum over all $m $ in $\M$ in the right hand side
of (\ref{oraclemulti}) at the price of the additional term $w_m$.
This result can be viewed as an oracle type property: indeed,
without knowing $(f-g)$, we know that the uniform separation rate of
the aggregated test is of the same order as the smallest uniform
separation rate in the collection of single tests, up to the factor
$w_m$. It will be used to prove that our multiple testing procedures
are adaptive over various classes of alternatives.

\smallskip

We focus here on two particular examples. The first example
involves a nested collection of linear subspaces of $
\mathbb{L}^2([0,1])$, as in model selection estimation approaches.
In the second example, we consider a collection of one dimensional
linear subspaces of $ \mathbb{L}^2([0,1])$, and our testing
procedure is hence related to a thresholding estimation approach.

\smallskip

\emph{[Multiple kernels case - Example 1]} Let $\X=[0,1]$ and $\nu$ be the Lebesgue measure on $[0,1]$. Let
$\{\p_0,\ \p_{(j,k)},\ j\in\mathbb{N},\ k\in\{0,\ldots,2^j-1\}\}$ be
the Haar basis of $\mathbb{L}^2([0,1])$ with
\begin{equation}\label{Haar}
\p_0(x)=\1_{[0,1]}(x) \quad \textrm{and} \quad
\p_{(j,k)}(x)=2^{j/2}\psi(2^jx-k),
\end{equation}
where $\psi(x)=\1_{[0,1/2)}(x)-\1_{[1/2,1)}(x)$. The collection of
linear subspaces $\{S_m,m\in\M\}$ is chosen as a collection of
nested subspaces generated by subsets of the Haar basis. More
precisely, we denote by $S_0$ the subspace of $\mathbb{L}^2([0,1])$
generated by $\p_0$, and we define $K_0(x,x')=\p_0(x)\p_0(x')$. We
also consider for $J\geq 1$ the subspaces $S_J$ generated by
$\{\p_\la, \la \in \{0\}\cup \La_J\}$ with $\La_J=\{(j,k),\
j\in\{0,\ldots, J-1\},\ k\in\{0,\ldots,2^j-1\}\}$, and
$K_J(x,x')=\sum_{\la\in\{0\}\cup\La_J} \p_\la(x)\p_\la(x')$. Let for
some $\bar{J}\geq 1$, $\M_{\bar{J}}=\{J, 0\leq J\leq \bar{J}\},$ and for every $J$ in $\M_{\bar{J}}$, $w_J=2 \left(\ln(J+1) + \ln (\pi/\sqrt{6})\right).$

Let $\Phi_{\alpha}^{(1)}$ be the test defined by
(\ref{testmult}) with the collection of kernels $\{K_J,J\in\M_{\bar{J}}\}$  and with $\{w_J, J\in\M_{\bar{J}}\}$. We
obtain from Theorem \ref{testmulti} that there exists
$C\pa{\alpha,\beta, \norm{f}_\infty, \norm{g}_\infty}>0$ such that $\P_{f,g}\pa{\Phi_{\alpha}^{(1)}=0} \leq \beta$
if
\begin{multline}\label{nested}
  \|f-g\|^2 \geq C\pa{\alpha,\beta, \norm{f}_\infty, \norm{g}_\infty} \inf_{J\in\M_{\bar{J}}}\Bigg\{
  \|(f-g)-\Pi_{S_J}(f-g)\|^2\\
 + (\ln (J+2)) \frac{2^{J/2}}{n}\Bigg\}.
  \end{multline}

\emph{[Multiple kernels case - Example 2]} Let $\X=[0,1]$ and $\nu$ be the Lebesgue measure on $[0,1]$. Let
$\{\p_0,\p_{(j,k)},j\in\mathbb{N},k\in\{0,\ldots,2^j-1\}\}$ still be
the Haar basis of $\mathbb{L}^2([0,1])$
 defined by (\ref{Haar}). Let for some $\tilde{J}\geq 1$,
 $$\La_{\tilde{J}}=\{(j,k),\ j\in\{0,\ldots, \tilde{J}-1\},\
k\in\{0,\ldots,2^j-1\}\}.$$ For any $\lambda$ in $\{0\}\cup
\Lambda_{\tilde{J}}$, we consider the subspace $\tilde{S}_\lambda$
of $\mathbb{L}^2([0,1])$ generated by $\p_\lambda$, and
$K_\lambda(x,x')=\p_\lambda(x)\p_\lambda(x')$. Let
$\Phi_{\alpha}^{(2)}$ be the test defined by (\ref{testmult}) with
the collection of kernels $\{K_\lambda,\lambda\in\{0\}\cup
\Lambda_{\tilde{J}}\}$, with $w_0=\ln(2)$, and $w_{(j,k)}=\ln(2^j)+2
\left(\ln(j+1) + \ln (\pi/\sqrt{3})\right)$ for
$j\in\N,k\in\{0,\ldots,2^j-1\}.$ We obtain from Theorem
\ref{testmulti} and Pythagoras' theorem that there is some constant
$C\pa{\alpha,\beta, \norm{f}_\infty, \norm{g}_\infty}>0$ such that
if there exists $\lambda$ in $\{0\}\cup \Lambda_{\tilde{J}}$ for
which
$$\|\Pi_{\tilde{S}_\lambda}(f-g)\|^2 \geq C\pa{\alpha,\beta, \norm{f}_\infty, \norm{g}_\infty} \frac{w_\lambda}{n},$$
then  $\P_{f,g}\pa{\Phi_{\alpha}^{(2)}=0 } \leq \beta$. If $\M_{\tilde{J}}=\{m,m\subset \{0\}\cup \Lambda_{\tilde{J}}\},$
the above condition is equivalent to saying that there exists $m$ in
$\M_{\tilde{J}}$ such that
$$\|\Pi_{S_m}(f-g)\|^2 \geq C\pa{\alpha,\beta, \norm{f}_\infty, \norm{g}_\infty} \frac{\sum_{\lambda\in m}
w_\lambda}{n},$$ where $S_m$ is generated by
$\{\p_\lambda, \lambda\in m\}$. Hence, there exists some constant
$C\pa{\alpha,\beta, \norm{f}_\infty,
\norm{g}_\infty}~>0$ such that $\P_{f,g}\pa{\Phi_{\alpha}^{(2)}=0} \leq \beta$ if
\begin{multline}\label{thresh}
\|f-g\|^2 \geq C\pa{\alpha,\beta, \norm{f}_\infty, \norm{g}_\infty}
\inf_{m\in \M_{\tilde{J}}}\Big\{ \|(f-g)-\Pi_{{S_m}}(f-g)\|^2 \\
+ \frac{\sum_{\lambda\in m} w_\lambda}{n}\Big\}.
\end{multline}

\subsubsection{Multiple testing procedures based on approximation kernels}\label{sectionmultikernel}

\begin{Th}\label{testmultapproxkern}
Let $\alpha,\beta$ be fixed levels in $(0,1)$, $\X=\R^d$ and let
 $\nu$ be the Lebesgue measure on $\R^d$.  Let $\{k_{m_1},m_1\in\M_1\}$ be a collection of approximation
kernels
 such that $\int_{\X} k_{m_1}^2(x)d\nu_x<+\infty $, $k_{m_1}(x)=k_{m_1}(-x)$, and a collection $\{h_{m_2},m_2\in \M_2\}$, where each $h_{m_2}$ is a vector of $d$ positive bandwidths $(h_{m_2,1},\ldots,h_{m_2,d})$.
 We set $\M=\M_1\times\M_2$, and for all $m=(m_1,m_2)$ in $\M$, $x=(x_1,\ldots,x_d)$, $x'=(x_1',\ldots,x_d')$ in
 $ \R^d$,
 $$K_m(x,x')= k_{m_1,h_{m_2}}(x-x')=\frac{1}{\prod_{i=1}^d
 h_{m_2,i}}k_{m_1}\pa{\frac{x_1-x_1'}{h_{m_2,1}},\ldots,\frac{x_d-x_d'}{h_{m_2,d}}}.$$
 Let $\Phi_{\alpha}$ be the
test defined by (\ref{testmult}) with $\{K_m,m\in\M\}$ and a
collection $\{w_m,m\in \M\}$ of positive numbers such that
$\sum_{m\in \M} e^{-w_m}\leq 1$.

Then $\Phi_{\alpha}$ is a level
$\alpha$ test. Moreover, there
 exists $\kappa>0$ such that
if
\begin{multline}\label{oraclenoyau}
\norm{f-g}^2\!\!\geq \!\!\inf_{(m_1,m_2)\in \M}
\left\{\norm{(f-g)-\!\! k_{m_1,h_{m_2}}\!\!*(f-g)}^2+\right.\nonumber\\
\left.\frac{4+2\sqrt{2}\kappa
(\ln(2/\alpha)+w_m)}{n\sqrt{\beta}}\sqrt{\frac{\norm{f+g}_\infty\norm{f+g}_1\norm{k_{m_1}}^2}{
\prod_{i=1}^dh_{m_2,i}}}\right\}+ \frac{8\norm{f+g}_\infty}{\beta n
},
\end{multline}
 then
$$\P_{f,g}\pa{\Phi_{\alpha}=0} \leq \beta.$$

\end{Th}

We focus here on two particular examples. The first example involves a collection of non necessarily integrable approximation kernels with a collection of bandwidths vectors whose components are the same in every direction. The second example involves a single integrable approximation kernel, but with a collection of bandwidths vectors whose components may differ according to every direction.

\smallskip

\emph{[Multiple kernels case - Example 3]} Let $\X= \R^d$ and $\nu$ be the Lebesgue measure on $\R^d$. We set
$\M_1=\N\setminus\{0\}$ and $\M_2=\N$. For $m_1$ in $\M_1$, let
$k_{m_1}$ be a kernel such that $\int k_{m_1}^2(x)d\nu_x<+\infty$
and $k_{m_1}(x)=k_{m_1}(-x)$, non necessarily integrable, whose
Fourier transform is defined when $k_{m_1}\in \LL^1(\R^d)\cap
\LL^2(\R^d)$ by $\widehat{k_{m_1}}(u)=\int_{\R^d} k_{m_1}(x)
e^{i\langle x,u\rangle}d\nu_x$ and is extended to  $k_{m_1}\in
\LL^2(\R^d)$ in the Plancherel sense. We assume that for every $m_1$
in $\M_1$, $\norm{\widehat{k_{m_1}}}_\infty<+\infty$, and
\begin{equation}\label{nonintkernel}
\mbox{Ess sup}_{u\in\R^d\setminus \{0\}}\frac{|1-
\widehat{k_{m_1}}(u)|}{\norm{u}_d^{m_1}}\leq C,
\end{equation}
for some $C>0$, where $\norm{u}_d$ denotes the euclidean norm of
$u$. Note that the sinc kernel, the spline type kernel and Pinsker's
kernel given in \cite{Sacha} for instance satisfy this condition
which can be viewed as an extension of the integrability condition
(see \cite{Sacha} p. 26-27 for more details). For
$m_2$ in $\M_2$, let $h_{m_2}=(2^{-m_2},\ldots,2^{-m_2})$
 and for $m=(m_1,m_2)$ in $\M=\M_1\times\M_2$, let
$$K_m(x,x')=k_{m_1,h_{m_2}}(x-x')=\frac{1}{2^{-dm_2}}k_{m_1}\pa{\frac{x_1-x_1'}{2^{-m_2}},\ldots,\frac{x_d-x_d'}{2^{-m_2}}}
.$$ We take $w_{(m_1,m_2)}= 2\left(\ln(m_1(m_2+1))+
\ln(\pi^2/6)\right)$, so $\sum_{m\in\mathcal{M}} e^{-w_m}\leq 1$.
Let $\Phi_{\alpha}^{(3)}$ be the test defined by
(\ref{testmult}) with the collection of kernels
$\{K_m,m\in\M\}$ and $\{w_m,m\in\M\}$. We obtain from Theorem
\ref{testmultapproxkern} that there exists $C(\alpha,\beta)>0$ such that $\P_{f,g}\pa{\Phi_{\alpha}^{(3)}=0}
\leq \beta$ if
\begin{multline}\label{oraclenoyausob}
\norm{f-g}^2 \geq  C(\alpha,\beta)\Bigg(\inf_{(m_1,m_2)\in \M}
\Bigg\{\norm{(f-g)-k_{m_1,h_{m_2}}\!\!*(f-g)}^2\\
+\frac{w_{(m_1,m_2)}}{n}\sqrt{\frac{\norm{f+g}_\infty\norm{f+g}_1\norm{k_{m_1}}^2}{2^{-dm_2}}}\Bigg\}+\frac{\norm{f+g}_\infty}{n}\Bigg).\end{multline}

\emph{[Multiple kernels case - Example 4]} Let $\X= \R^d$ and $\nu$ be the Lebesgue measure on $\R^d$. Let
$\M_1=\{1\}$ and $\M_2=\N^d$. For $x=(x_1,\ldots,x_d)$ in $\R^d$,
let $k_1(x)=\prod_{i=1}^{d} k_{1,i}(x_i)$ where the $k_{1,i}$'s are
real valued kernels such that $k_{1,i}\in \LL^1(\R)\cap \LL^2(\R)$,
$k_{1,i}(x_i)=k_{1,i}(-x_i)$, and $\int_{\R} k_{1,i}(x_i)dx_i=1$.
For $m_2=(m_{2,1},\ldots,m_{2,d})$ in $\M_2$,
$h_{m_2,i}=2^{-m_{2,i}}$ and for $m=(m_1,m_2)$ in
$\M=\M_1\times\M_2$,
$$K_m(x,x')=k_{m_1,h_{m_2}}(x-x')=\prod_{i=1}^{d}
\frac{1}{h_{m_2,i}}k_{1,i}\left(\frac{x_i-x_i'}{h_{m_2,i}}\right).$$
We also set
 $w_{(1,m_2)}=2\sum_{i=1}^d \left(\ln
(m_{2,i}+1)+\ln(\pi/\sqrt{6})\right)$, so that\\
$\sum_{m\in\M_1\times\M_2} e^{-w_m}=1$. Let $\Phi_{\alpha}^{(4)}$ be
the  test defined by (\ref{testmult}) with the collections
$\{K_m,m\in\M\}$  and $\{w_m,m\in\M\}$. We deduce from
Theorem~\ref{testmultapproxkern} that there exists $C(\alpha,\beta)>0$ such that $\P_{f,g}\pa{\Phi_{\alpha}^{(4)}=0} \leq \beta$
if
\begin{multline}\label{oraclenoyaunik}
\norm{f-g}^2\geq
 C(\alpha,\beta)\Bigg(\inf_{m_2\in \M_2}
\Bigg\{\norm{(f-g)-k_{1,h_{m_2}}*(f-g)}^2\\
+\frac{w_{(1,m_2)}}{n}\sqrt{\frac{\norm{f+g}_\infty\norm{f+g}_1\norm{k_{1}}^2}{
\prod_{i=1}^d h_{m_2,i}}}\Bigg\} +\frac{\norm{f+g}_\infty}{n}\Bigg).
\end{multline}

\subsection{Uniform separation rates over various classes of alternatives}

We here evaluate the uniform separation rates, defined by
(\ref{defrate}), of the
multiple testing procedures introduced above over several classes of
alternatives based on Besov and weak Besov bodies when
$\X=[0,1]$, or Sobolev and anisotropic Besov-Nikol'skii balls
when $\X=\R^d$.

\subsubsection{Uniform separation rates for Besov and weak Besov bodies}\label{Besov}
In this section, we adapt to the present setting the results that we obtained in \cite{MBP}.

Given $\alpha$ in $(0,1)$, let $\Phi_{\alpha/2}^{(1)}$ and $\Phi_{\alpha/2}^{(2)}$ be the tests defined
in \emph{[Multiple kernels case - Example 1]} and \emph{[Multiple kernels case - Example 2]} (with $\alpha$ replaced by $\alpha/2$), and let $\Psi_\alpha=\max(\Phi_{\alpha/2}^{(1)},\Phi_{\alpha/2}^{(2)}).$\\
Recall that these tests are
constructed from the Haar basis
$\{\p_0,\p_{(j,k)},j\in\mathbb{N},k\in\{0,\ldots,2^j-1\}\}$ of
$\mathbb{L}^2([0,1])$
 defined by (\ref{Haar}).
We define for $\delta>0$, $R>0$ the Besov body $\mathcal{B}_{2,\infty}^\delta(R)$ as follows:
\begin{multline*}
\mathcal{B}_{2,\infty}^\delta(R)=\Bigg\{ s=\alpha_0\p_0+\sum_{j\in\mathbb{N}}\sum_{k=0}^{2^j-1} \alpha_{(j,k)} \p_{(j,k)}\ \Big/\ \alpha_0^2\leq R^2,\ \forall j
\in \mathbb{N},\\ \sum_{k=0}^{2^j-1} \alpha_{(j,k)}^2\leq R^2 2^{-2j\delta}\Bigg\}.
\end{multline*}
We also consider the weak Besov body given for $\gamma >0$, $R'>0$ by
\begin{multline*}
\mathcal{W}_{\gamma}(R')=\Bigg\{ s=\alpha_0\p_0+\sum_{j \in \mathbb{N}}\sum_{k=0}^{2^j-1}\alpha_{(j,k)}
\p_{(j,k)}\ \Big/\\
 \forall t>0,\ \alpha_0^2\1_{\alpha_0^2\leq t}+\sum_{j
\in \mathbb{N}} \sum_{k=0}^{2^j-1} \alpha_{(j,k)}^2\1_{\alpha_{(j,k)}^2\leq t}\leq R'^2  t^{\frac{2\gamma}{1+2\gamma}}\Bigg\}.
\end{multline*}

\begin{Cor}\label{vitesses}
Assume that $\ln\ln n\geq 1$, $2^{\bar{J}} \geq n^2$, and
$\tilde{J}=+\infty$. Then, for any $\delta>0$,  $\gamma >0$, $R,R',
R''>0$, if \begin{multline*}
\mathcal{B}_{\delta, \gamma,\infty}(R,R',R'')=\big\{(f,g)\ \big/ \
(f-g)\in
\mathcal{B}_{2,\infty}^\delta(R)\cap \mathcal{W}_{\gamma}(R'),\\
\max(\norm{f}_\infty,\norm{g}_\infty)\leq R''\big\},
\end{multline*}
$\rho(\Psi_\alpha, \mathcal{B}_{\delta, \gamma,\infty}(R,R',R'') ,\beta) $, defined by (\ref{defrate}), is upper bounded by \\
$(i)$ $C(\delta, \gamma, R,R',R'', \alpha, \beta)
\left(\frac{\ln\ln n}{n}\right)^{\frac{2\delta}{4\delta+1}}$ if $\delta\geq \gamma/2$,\\
$(ii)$ $C(\delta, \gamma, R,R',R'', \alpha, \beta)
\left(\frac{\ln n}{n}\right)^{\frac{\gamma}{2\gamma+1}}$ if $\delta < \gamma/2$.
\end{Cor}

\emph{Comments.}

1. Lower bounds for the minimax separation rates over
 $\mathcal{B}_{\delta, \gamma,\infty}(R,R',R'')$ are also available, proving that the test $\Psi_\alpha$ is adaptive in the minimax sense over
$\mathcal{B}_{\delta, \gamma,\infty}(R,R',R'')$,
up to a $\ln\ln n$ factor if $ \delta \geq \max{}(\gamma/2,\gamma/(1+2\gamma))$ and exactly if $\delta<\gamma/2$ and $\gamma>1/2$.
In the other cases, the  exact rate is unknown.

2. Let us mention here that our classes of alternatives are not
defined in the same way as in \cite{ButuceaTribouley} in the
classical two-sample problem for i.i.d. samples, since the classes of
alternatives $(f,g)$ of \cite{ButuceaTribouley} are such that $f$
and $g$ both belong to a Besov ball. Here the smoothness condition
is only required on the difference $(f-g)$.  In
particular, the functions $f$ and $g$ might be very irregular but as
long as their difference is smooth, the probability of second kind error of the test will be controlled.

\subsubsection{Uniform separation rates for Sobolev and anisotropic Nikol'skii-Besov balls}

Let $\Phi_{\alpha}^{(3)}$ be defined as in \emph{[Multiple kernels case - Example 3]}, and let us introduce for $\delta>0$ the Sobolev ball
$\mathcal{S}_d^\delta(R)$ defined by
$$\mathcal{S}_d^\delta(R)=\Bigg\{s:\R^d\to\R\ \Big/
s\in\mathbb{L}^1(\R^d)\cap\mathbb{L}^2(\R^d),\ \int_{\R^d}
\norm{u}_d^{2\delta}|\hat{s}(u)|^2du\leq (2\pi)^d R^2\Bigg\},$$
where $\norm{u}_d$ denotes the euclidean norm of $u$ and $\hat{s}$
denotes the Fourier transform of $s$: $\hat{s}(u)=\int_{\R^d}
s(x)e^{i\langle x,u\rangle} dx.$

\begin{Cor}\label{vitessessobolev}
Assume that $\ln \ln n\geq 1$. For any $\delta, R, R', R''>0$, if
\begin{multline*}
\mathcal{S}_d^{\delta}(R,R',R'')=\{(f,g)\ \big/
(f-g)\in\mathcal{S}_{d}^{\delta}(R),\
\max(\norm{f}_1,\norm{g}_1)\leq R',\\
\max(\norm{f}_\infty,\norm{g}_\infty)\leq R''\},
\end{multline*}
then
$$\rho(\Phi_{\alpha}^{(3)},\mathcal{S}_{d}^{\delta}(R,R',R''),\beta)\leq
C(\delta,\alpha,\beta,R,R',R'',d)\left(\frac{\ln\ln
n}{n}\right)^{\frac{2\delta}{d+4\delta}}.$$
\end{Cor}

\emph{Comments.} From \cite{RigolletTsybakov}, we know that, in the
density model, the minimax adaptive estimation rate over
$\mathcal{S}_{d}^{\delta}(R)$   is of order
$n^{-\frac{\delta}{d+2\delta}}$ when $\delta>d/2$. Rigollet and
Tsybakov construct some aggregated  density estimators, based on
Pinsker's kernel, that achieve this rate with exact constants. In
the same way, the test $\Phi_{\alpha}^{(3)}$ consists in an
aggregation of some tests based on a collection of kernels, that may
be for instance a collection of Pinsker's kernels. It achieves over
$\mathcal{S}_{d}^{\delta}(R,R',R'')$ a uniform separation rate of
order $n^{-\frac{2\delta}{d+4\delta}}$ up to a $\ln\ln n$ factor.
This rate is now known to be the optimal adaptive minimax rate of
testing when $d=1$ in several models (see \cite{Spokoiny96} in a
Gaussian model or \cite{Ingster2000} in the density model  for
instance). From the results of \cite{HorowitzSpokoiny2001}, we can
conjecture that our rates are also optimal when $d>1$.

\smallskip

Let $\Phi_{\alpha}^{(4)}$ be the test defined in \emph{[Multiple kernels case - Example 4]}. Let $\Delta=(\Delta_1,\ldots,\Delta_d)$, where for every $i=1\ldots d$, $\Delta_i$ is a positive integer. Assume furthermore that  $\int_{\R}
|k_{1,i}(x_i)| |x_i|^{\Delta_i }dx_i<+\infty$, and $\int_{\R}
k_{1,i}(x_i)x_i^jdx_i=0$ for every $i=1\ldots d$ and $j=1\ldots \Delta_i$.

For $\delta=(\delta_1,\ldots,\delta_d)\in \prod_{i=1}^d(0,\Delta_i]$ and $R>0$, we consider the anisotropic Nikol'skii-Besov ball
$\mathcal{N}_{2,d}^{\delta}(R)$ defined by:
\begin{multline*}
 \mathcal{N}_{2,d}^\delta(R)=\Big\{s:\R^d\to \R\ \big/\  s \textrm{ has continuous partial derivatives } D_i^{\lfloor \delta_i\rfloor}\\
 \textrm{ of order }
 \lfloor \delta_i\rfloor \textrm{ w.r.t } u_i,
\textrm{ and }  \forall i=1\ldots
d,\ u_1,\ldots,u_d,v\in\R,\\
\norm{D_i^{\lfloor \delta_i \rfloor
}s(u_1,\ldots,u_i+v,\ldots,u_d)-D_i^{\lfloor \delta_i
\rfloor}s(u_1,\ldots,u_d) }_2
\leq R|v|^{\delta_i-\lfloor \delta_i \rfloor} \Big\}.
\end{multline*}

\begin{Cor}\label{vitessesanisotrop}
Assume that $\ln \ln n\geq 1$. For any $\delta=(\delta_1,\ldots,\delta_d)$ in $\prod_{i=1}^d (0,\Delta_i]$
and $R, R',R''>0$, if
\begin{multline*}
\mathcal{N}_{2,d}^{\delta}(R,R',R'')=\{(f,g)\ \big/\
(f-g)\in\mathcal{N}_{2,d}^{\delta}(R),\
\max(\norm{f}_1,\norm{g}_1)\leq
R',\\\max(\norm{f}_\infty,\norm{g}_\infty)\leq R''\},
\end{multline*}
then, for
${1}/{\bar{\delta}}=\sum_{i=1}^d {1}/{\delta_i}$,
$$\rho(\Phi_{\alpha}^{(4)},\mathcal{N}_{2,d}^{\delta}(R,R',R''),\beta)\leq
C(\delta,\alpha,\beta,R,R',R'',d) \left(\frac{\ln \ln
n}{n}\right)^{\frac{2\bar{\delta}}{1+4\bar{\delta}}}.$$
\end{Cor}

\emph{Comments.} When $d=1$, from \cite{Ingster2000}, we know that
in the density model, the adaptive minimax rate of testing  over
 a Nikol'skii class with smoothness parameter $\delta$ is of order $(\ln \ln
n/n)^{2\delta/(1+4\delta)}$. We find here an upper bound similar to
this univariate rate, but where $\delta$ is replaced by
$\bar{\delta}$. Such results were obtained in a multivariate density
estimation context in \cite{GoldLep} where the adaptive minimax
estimation rates over the anisotropic Nikol'skii classes are proved
to be of order $n^{-\bar{\delta}/(1+2\bar{\delta})}$, and where
adaptive kernel density estimators are proposed. Moreover, the
minimax rates of testing obtained recently in
\cite{IngsterStepanova2011} over anisotropic periodic Sobolev balls,
but in the Gaussian white noise model, are of the same order as the
upper bounds obtained here.

\section{Proofs}\label{preuves}

\subsection{Proof of Proposition \ref{marques}}
All along the proof, $\int$ denotes $\int_{\X}$. Recalling that the marked point processes are characterized by their Laplace functional (see \cite{daleyVerejones} for instance), we first aim at computing
$\E\left[\exp\left(\int h dN\right)\right]$ for any bounded measurable
function $h$ on $\X$.
Since $N^1$ and $N^{-1}$ are independent,
$$\E\left[\exp\left(\int h dN\right)\right]=\E\left[\exp\left(\int h dN^1\right)\right]\E\left[\exp\left(\int h dN^{-1}\right)\right].$$
The Laplace functional of $N^1$ is given by
$$\E\left[\exp\left(\int h dN^1\right)\right]=\exp\left(\int
\left(e^h-1\right) f d\mu\right),$$ and the Laplace functional of
$N^{-1}$ has the same form, replacing $f$ by $g$, so
$$\E\left[\exp\left(\int h dN\right)\right]=\exp\left(\int \left(e^h-1\right) \left(f+g\right) d\mu\right),$$
which is the Laplace functional of a Poisson process with intensity $(f+g)$ w.r.t. $\mu$. Therefore, $N$ is a Poisson process with intensity $(f+g)$ w.r.t. $\mu$.

In order to prove (\ref{loieps0}), we then give an explicit expression of the function:
$$ t=(t_x)_{x\in N} \mapsto  \Phi(t,N)=\E\left[\exp\left(\sum_{x\in N} t_x \e^0_x \right) \Bigg|N\right],$$
which characterizes the distribution of $(\e^0_x)_{x\in
N}$ conditionally on $N$.\\
Let $\lambda$ be a bounded measurable function defined on $\X$, and let
$$\E_{\lambda}=\E\cro{\exp\left(\int \lambda dN \right)\exp\left(\sum_{x\in N} t_x \e^0_x\right)}.$$
By definition of $(\e^0_x)_{x\in
N}$ and by independency of $N^1$ and $N^{-1}$, we have that
\begin{eqnarray*}
 \E_{\lambda}&=&\E\cro{\exp\pa{\int (\lambda(x)+t_x) dN^1_x} \exp\pa{\int (\lambda(x)-t_x) dN^{-1}_x} }\\
&=& \E\cro{\exp\pa{\int (\lambda(x)+t_x) dN^1_x} }
\E\cro{\exp\pa{\int (\lambda(x)-t_x) dN^{-1}_x} }\\
&=&\exp \int \cro{( e^{\lambda(x)+t_x}-1 )f(x)+
 ( e^{\lambda(x)-t_x}-1 )g(x)} d\mu_x.
 \end{eqnarray*}
Then, for $h(x)= \lambda(x)+\ln\pa{ \frac{e^{t_x} f(x) + e^{-t_x} g(x)}{(f+g)(x)}}$,
$$\E_{\lambda}=\exp \int (e^{h(x)}-1) (f+g)(x) d\mu_x=\E\cro{\exp\left(\int h dN \right)}.$$
Hence, for every bounded measurable function $\lambda$ defined on
$\X$,
\begin{multline*}
\E\cro{\exp\left(\int \lambda dN \right)\exp\left(\sum_{x\in N} t_x
\e^0_x\right)}=\\
 \E \Bigg[\exp\left(\int \lambda dN\right)
\prod_{x\in N} \Bigg(e^{t_x} \frac{f(x)}{(f+g)(x)}
+ e^{-t_x} \frac{g(x)}{(f+g)(x)}\Bigg)\Bigg].
\end{multline*}
Since the marked point processes are characterized by their Laplace
functional, this implies that
$$\Phi(t,N)=\E\cro{\exp\left(\sum_{x\in N} t_x \e^0_x\right) \Bigg|N}  = \prod_{x\in N} \pa{e^{t_x} \frac{f(x)}{(f+g)(x)} +
e^{-t_x} \frac{g(x)}{(f+g)(x)}},$$
which concludes the proof.

\subsection{Proof of Proposition
\ref{propsinglerror}}\label{proofprop2}

%\begin{proof}

Let us prove the first part of Proposition~\ref{propsinglerror}.
Recall that $q_{K,1-\beta/2}^\alpha$ denotes the $1-\beta/2$ quantile
of $q^{(N)}_{K,1-\alpha}$, which is the $(1-\alpha)$ quantile of
$\hat{T}^{\e}_K$ conditionally on $N$. We here want to find a condition on $\hat T_K$, or more precisely on
$\mathcal{E}_K=\E_{f,g}[\hat T_K]$, ensuring that
$$\P_{f,g}(\hat{T}_K\leq q_{K,1-\beta/2}^\alpha)\leq \beta/2.$$
From
Markov's inequality, we have that for any $x>0$,
$$\P_{f,g}\left(\left|-\hat T_K+\mathcal{E}_K\right|\geq x\right)\leq \frac{\Var(\hat T_K)}{x^2}.$$
Let us compute $\Var(\hat T_K)=\E_{f,g}[\hat T_K^2]-\mathcal{E}_K^2$.
Let $\X^{[3]}$ and $\X^{[4]}$ be the sets
$\{(x,y,u)\in\X^3,\ x,y,u \textrm{ all different}\}$ and
$\{(x,y,u,v)\in\X^4,\ x,y,u,v \textrm{ all different}\}$ respectively.
Since
$$\E_{f,g}[\hat T_K^2]=\E_{f,g}\left[\E\left[\left(\int_{\X^{[2]}} K(x,x')\e_x^0\e_{x'}^0 dN_xdN_{x'}\right)^2\Bigg|N\right]\right],$$
by using (\ref{loieps0}),
\begin{eqnarray*}
\E_{f,g}[\hat T_K^2]&=&\E_{f,g}\Bigg[\int_{\X^{[4]}}\! K(x,y)K(u,v) \frac{f-g}{f+g}(x) \frac{f-g}{f+g}(y)\\
&&\frac{f-g}{f+g}(u)\frac{f-g}{f+g}(v)dN_xdN_ydN_udN_v\Bigg]\\
&&+4\E_{f,g}\left[\int_{\X^{[3]}}  K(x,y)K(x,u) \frac{f-g}{f+g}(y)
\frac{f-g}{f+g}(u)dN_xdN_ydN_u\right]\\
&& +2\E_{f,g}\left[\int_{\X^{[2]}} K^2(x,y)dN_xdN_y\right].
\end{eqnarray*}
Now, from Lemma 5.4 III in \cite{daleyVerejones} on factorial moments measures applied to Poisson
processes, we deduce that
\begin{eqnarray*}
\E_{f,g}[\hat T_K^2]&=&\int_{\X^{4}}\Bigg( K(x,y)K(u,v)(f-g)(x)(f-g)(y)\\
&&(f-g)(u)(f-g)(v)\Bigg)  d\mu_xd\mu_yd\mu_ud\mu_v\\
&&+4\int_{\X^{3}}  K(x,y)K(x,u)(f+g)(x)(f-g)(y)(f-g)(u) d\mu_xd\mu_yd\mu_u\\
&&+2\int_{\X^{2}} K^2(x,y)(f+g)(x)(f+g)(y) d\mu_xd\mu_y\\
\end{eqnarray*}
Note that the three above
integrals are finite, thanks to Assumptions \ref{Hypo1}, \ref{Hypo2} et \ref{Hypo3}. We
finally obtain that
$\E_{f,g}[\hat T_K^2]=\mathcal{E}_K^2+4 n^3 A_K+2 n^2 B_K,$ and for $x>0$,
$$\P_{f,g}\left(\left|-\hat T_K+\mathcal{E}_K\right|\geq x\right)\leq \frac{4n^3 A_K+2 n^2B_K}{x^2}.$$
Taking $x=2n \sqrt{(2nA_K+B_K)/\beta}$ in the above inequality leads to
\begin{equation}\label{ineMarkov}
\P_{f,g}\left(\left|-\hat T_K+\mathcal{E}_K\right|\geq 2 n\sqrt{\frac{2n A_K+B_K}{\beta}}\right)\leq \frac{\beta}{2}.
\end{equation}
Therefore, if
$\mathcal{E}_K>2n \sqrt{\frac{2n A_K+B_K}{\beta}}+q_{K,1-\beta/2}^\alpha$,
then $\P_{f,g}(\hat{T}_K\leq q_{K,1-\beta/2}^\alpha)\leq \beta/2, $ so
 $\P_{f,g}(\Phi_{K,\alpha}=0)\leq \beta.$

\smallskip

Let us now give a sharp upper bound for $q_{K,1-\beta/2}^\alpha$. Reasoning conditionally on $N$, we recognize in $\hat{T}^{\e}_K$ a homogeneous Rademacher chaos, as defined by de la Pe$\tilde{\mbox{n}}$a and Gin\'e \cite{PenaGine},  of the form $X=\sum_{i\neq i'} x_{i,i'} \e_i\e_{i'},$
where the $x_{i,i'}$'s are some real deterministic numbers and $(\e_i)_{i\in\N}$ is a sequence of i.i.d. Rademacher variables. Corollary 3.2.6 of  \cite{PenaGine} states that there exists some absolute constant $\kappa>0$ such that if $\sigma^2=\E[X^2]=\sum_{i\neq i'} x_{i,i'}^2,$ then
$$\E\left[\exp\left(|X|/(\kappa \sigma)\right)\right] \leq 2.$$
Hence by Markov's inequality,
$$\P(|X|\geq \kappa \sigma \ln(2/\alpha))\leq \alpha.$$
Note that one could find more precise constants with the results of \cite{Latala}.\\
Applying this result to $\hat{T}^{\e}_K$ with $\sigma^2=\sum_{x\neq x'\in N } K^2(x,x')$
leads to
\begin{equation*}
q^{(N)}_{K,1-\alpha}\leq \kappa \ln(2/\alpha)
\sqrt{\int_{\X^{[2]}} K^2(x,y)dN_x dN_{y}}.
\end{equation*}
Hence $q_{K,1-\beta/2}^\alpha$ is upper bounded by the $(1-\beta/2)$
quantile of\\
 $\kappa \ln(2/\alpha) \sqrt{\int_{\X^{[2]}}
K^2(x,y)dN_x dN_{y}}$.\\
Using Markov's inequality again and  Lemma~5.4~III  in \cite{daleyVerejones}, we obtain that $$\P_{f,g}\left(\int_{\X^{[2]}} K^2(x,y)dN_x dN_{y}\geq \frac{2 n^2 B_K}{\beta}\right)\leq \frac{\beta}{2},$$
and
\begin{equation*}
q_{K,1-\beta/2}^\alpha \leq \kappa \ln\left(2/\alpha\right) n \sqrt{\frac{2 B_K}{\beta}}.
\end{equation*}

\subsection{Proof of Theorem \ref{singlerror}}

First notice that for every $r>0$, and every kernel function $K$ satisfying Assumption \ref{Hypo3},
$$\mathcal{E}_K=\frac{n^2r}{2}\left(\norm{f-g}^2+r^{-2}\norm{K\cro{f-g}}^2-\norm{(f-g)-r^{-1}
K\cro{f-g}}^2\right).$$
With the notations of Proposition \ref{propsinglerror}, let $C_K$ be any upper bound for $B_K$. Since $A_K\leq  \norm{K\cro{f-g}}^2 \norm{f+g}_\infty$, from Proposition \ref{propsinglerror}, we deduce that $\P_{f,g}\left(\Phi_{K,\alpha}=0\right)\leq \beta$ if
\begin{multline*}
\norm{f-g}^2+r^{-2}\norm{K\cro{f-g}}^2-\norm{(f-g)-r^{-1} K\cro{f-g}}^2\\
\geq
4\sqrt{\frac{2\norm{f+g}_\infty}{n\beta}}\frac{\norm{K\cro{f-g}}}{r}+\frac{2}{nr\sqrt{\beta}}\left(2+\kappa\sqrt{2}\ln\left(\frac{2}{\alpha}\right)\right)
\sqrt{C_K},
\end{multline*}
By using the elementary inequality $2ab\leq a^2+b^2$ with
$a=\norm{K\cro{f-g}}/r$ and
$b=2\sqrt{2}\sqrt{\norm{f+g}_\infty/(n\beta)}$ in the right hand
side of the above condition, this condition can be replaced by:
\begin{multline*}
\norm{f-g}^2\geq \norm{(f-g)-r^{-1}
K\cro{f-g}}^2+\frac{8\norm{f+g}_\infty}{n\beta}\\
+\frac{2}{nr\sqrt{\beta}}\left(2+\kappa\sqrt{2}\ln\left(\frac{2}{\alpha}\right)\right)
\sqrt{C_K}. \end{multline*}
We can even add an infimum over $r$ in the right hand side of the condition, since $r$ can be arbitrarily chosen. Let us now justify our choices for $C_K$.

\smallskip

\emph{[Projection kernel case]} We consider  an orthonormal basis $\{\p_\la, \la \in \Lambda\}$ of a subspace $S$ of
 $\LL^{2}(\X,d\nu)$ and $K(x,x')=\sum_{\la\in\La}\p_\la(x)\p_\la(x').$
When the dimension of $S$ is finite, equal to $D$,
\begin{eqnarray*}
B_K&\leq& \norm{f+g}_\infty^2\int_\X \left(\sum_{\la\in\La}\p_\la(x)\p_\la(x')\right)^2d\nu_xd\nu_{x'}\\
&\leq&\norm{f+g}_\infty^2 D.
\end{eqnarray*}
When the dimension of $S$ is infinite,
\begin{eqnarray*}
B_K&= &   \int_{\X^2} \left(\sum_{\la\in\La}\p_\la(x)\p_\la(x')\right)^2(f+g)(x)(f+g)(x')  d\nu_xd\nu_{x'}\\
&\leq& \norm{f+g}_\infty  \int_{\X^2} \left(\sum_{\la\in\La}\p_\la(x)\p_\la(x')\right)^2(f+g)(x')  d\nu_xd\nu_{x'}\\
&\leq& \norm{f+g}_\infty  \int_{\X^2} \left(\sum_{\la, \la'
\in\La}\p_\la(x)\p_\la(x')\p_{\la'}(x)\p_{\la'}(x')\right)
(f+g)(x')  d\nu_xd\nu_{x'}\\
&\leq& \norm{f+g}_\infty  \sum_{\la ,\la' \in\La}  \int_\X \p_\la(x)\p_{\la'}(x)d\nu_x
 \int_\X  \p_\la(x') \p_{\la'}(x') (f+g)(x')  d\nu_{x'},
\end{eqnarray*}
where we have used the assumption (\ref{AA2}) to invert the sum and the
integral. Hence we have, by orthogonality, and since by assumption (\ref{AA1}), $\sum_{\la
\in\La}  \p^2_\la(x) \leq D$,
\begin{eqnarray*}
B_K&\leq& \norm{f+g}_\infty  \sum_{\la  \in\La}  \int_\X  \p^2_\la(x') (f+g)(x')  d\nu_{x'} \\
&\leq& \norm{f+g}_\infty  \norm{f+g}_1 D.
\end{eqnarray*}

\emph{[Approximation kernel case]} Assume now that $\X=\R^d$ and introduce an approximation kernel such
that $\int k^2(x) d\nu_x <+\infty$ and $k(-x)=k(x)$,
$h=(h_1,\ldots,h_d)$, with $h_i>0$ for every $i$, and $K(x,x')=k_h
(x-x')$, with $k_h(x_1,\ldots,x_d)=\frac{1}{\prod_{i=1}^d
h_i}k\pa{\frac{x_1}{h_1},\ldots,\frac{x_d}{h_d}}$. In this case,
\begin{eqnarray*}
B_K&=&\int_\X k_h^2(x-x')(f+g)(x)(f+g)(x')d\nu_xd\nu_{x'}\\
&\leq & \norm{f+g}_\infty \int_\X k_h^2(x-x')(f+g)(x)d\nu_xd\nu_{x'},\\
&\leq &  \frac{ \norm{f+g}_\infty \norm{f+g}_{1} \norm{k}^2
}{\prod_{i=1}^d h_i} .
\end{eqnarray*}
This ends the proof of Theorem \ref{singlerror}.

\subsection{Proof of Theorem \ref{vitesseparametrique}}

We first recall that when $K$ is chosen as in the \emph{[Reproducing kernel case]}, under the assumptions of Theorem \ref{vitesseparametrique}, $\mathcal{E}_K = n^2  \left\|m_f-m_g\right\|_{\mathcal{H}_K}^2$ (see Section \ref{introker}).

Since $A_K=  \int_{\X} \left\langle \int_{\X}  \theta(x)(f-g)(x) d\nu_x ,  \theta(y)\right\rangle_{\mathcal{H}_K}^2(f+g)(y) d\nu_y$, by the Cauchy-Schwarz inequality for the norm $\|.\|_{\mathcal{H}_K}$ in the RKHS, we obtain:
$$A_K\leq  \int_{\X} \left\|\int_{\X}  \theta(x)(f-g)(x) d\nu_x\right\|_{\mathcal{H}_K}^2 \norm{\theta(y)}_{\mathcal{H}_K}^2 (f+g)(y) d\nu_y.$$
Now, since for every $y$ in $\X$, $\|\theta(y)\|_{\mathcal{H}_K}^2=K(y,y)=\kappa_0$,
\begin{eqnarray*}
A_K&\leq &\kappa_0 \left\|\int_{\X}  \theta(x)(f-g)(x) d\nu_x\right\|_{\mathcal{H}_K}^2 \norm{f+g}_1\\
&\leq & \kappa_0 \left\|m_f-m_g\right\|_{\mathcal{H}_K}^2 \norm{f+g}_1\\
\end{eqnarray*}
This leads to
\begin{eqnarray*}
 2n\sqrt{\frac{2n A_K}{\beta}} &\leq&  2n\sqrt{\frac{2\kappa_0n\norm{f+g}_1}{\beta}}\left\|m_f-m_g\right\|_{\mathcal{H}_K} \\
&\leq& \frac{n^2}{2}\left\|m_f-m_g\right\|_{\mathcal{H}_K}^2 + 4 \frac{\kappa_0n\norm{f+g}_1}{\beta}.
\end{eqnarray*}
Finally, noting that $B_K\leq \kappa_0^2 \norm{f+g}_1^2$ and that by assumption $\norm{f+g}_1=2$, we obtain the desired result from Proposition \ref{propsinglerror} and obvious calculations.

\subsection{Proof of Proposition \ref{MC1}}
First let us rewrite here a result due to Romano and Wolf \cite{RW}.
\begin{Lemma}\label{LemmeRW}
\label{RomW}
Let $Y_0,...,Y_B$ be $B+1$ exchangeable variables then for all $u\in [0,1]$,
$$\P\pa{\frac{1}{B+1}\pa{1+\sum_{i=1}^B {\bf 1}_{Y_i\geq Y_0}}\leq u}\leq u.$$
\end{Lemma}
Assume that $(H_0)$ is satisfied. Conditionally on $N$, the observed statistic $\hat{T}_K^{\varepsilon^0}:=\hat{T}_K$ has the same distribution and is independent of the $\hat{T}_K^{\varepsilon^b}$'s for $b=1,...,B$. Therefore the variables $\hat{T}_K^{\varepsilon^b}$'s for $b=0,...,B$ are exchangeable variables given $N$.
Hence applying Lemma \ref{LemmeRW}, we obtain:
\begin{eqnarray*}
 \P\pa{ \hat{\Phi}_{\alpha}^K =1 \Big|N} &=& \P\pa{\hat{T}_K>\hat{q}^{(N)}_{K,1-\alpha}\Big|N}\\
&=& \P\pa{\sum_{b=1}^B {\bf 1}_{\hat{T}_K^{\varepsilon^b}\geq \hat{T}_K^{\varepsilon^0}} \leq \lfloor B \alpha \rfloor\Bigg|N }\\
&=& \P\pa{\frac{1}{B+1}\pa{1+\sum_{b=1}^B {\bf 1}_{\hat{T}_K^{\varepsilon^b}\geq \hat{T}_K^{\varepsilon^0}}} \leq \frac{\lfloor B \alpha \rfloor+1}{B+1}\Bigg|N }\\
&\leq & \frac{\lfloor B \alpha \rfloor+1}{B+1}.
\end{eqnarray*}

\subsection{Proof of Proposition \ref{puisMC}}
Let $t=q^{\alpha_B}_{K,1-\beta_B/2}$. By definition of $\hat{q}^{(N)}_{K,1-\alpha}$,
$$ \P_{f,g}( \hat{q}^{(N)}_{K,1-\alpha} >t) = \P_{f,g}(F_{K,B}(t) <1-\alpha) = \P_{f,g}\pa{\sum_{b=1}^B {\bf 1}_{ \hat{T}_K^{\varepsilon^b}\leq t}< B(1-\alpha)}.$$
We have
\begin{multline*}
\P_{f,g}\pa{\sum_{b=1}^B {\bf 1}_{ \hat{T}_K^{\varepsilon^b}\leq t}< B(1-\alpha),
 \ F_K(t)\geq 1-\alpha_B }\\
\leq  \P_{f,g}\pa{\sum_{b=1}^B \pa{{\bf 1}_{ \hat{T}_K^{\varepsilon^b}\leq t}-F_K(t)}< B(1-\alpha)-B (1-\alpha_B)}.
\end{multline*}
So we can decompose as follows:
\begin{multline*}
 \P_{f,g}\left( \hat{q}^{(N)}_{K,1-\alpha} >t\right)
  \leq \P_{f,g}(F_K(t)<1-\alpha_B)\\ +
 \P_{f,g}\pa{\sum_{b=1}^B \pa{{\bf 1}_{ \hat{T}_K^{\varepsilon^b}\leq t}-F_K(t)}< - B\sqrt{\frac{\ln B}{2B}}  }
\end{multline*}
By Hoeffding's inequality applied to the second probability given $N$, we obtain:
$$ \P_{f,g}( \hat{q}^{(N)}_{K,1-\alpha} >t) \leq \P_{f,g}(F_K(t)<1-\alpha_B) + \frac{1}{B}.$$
But by definition of $t$, this becomes
$$ \P_{f,g}( \hat{q}^{(N)}_{K,1-\alpha} >t) \leq \frac{\beta}{2}.$$
Let us now control the probability of second kind error of the test $\hat{\Phi}_{K,\alpha}$.
\begin{eqnarray*}
 \P_{f,g}( \hat T_K \leq  \hat{q}^{(N)}_{K,1-\alpha} )&\leq& \P_{f,g}( \hat T_K \leq \hat{q}^{(N)}_{K,1-\alpha}, \ \hat{q}^{(N)}_{K,1-\alpha}\leq t)+
\P_{f,g}(\hat{q}^{(N)}_{K,1-\alpha} >t)\\
&\leq&  \P_{f,g}( \hat T_K \leq t)+  \beta/2.
\end{eqnarray*}
We deduce from (\ref{ineMarkov}) that if
$$\mathcal{E}_K>2n \sqrt{\frac{2nA_K+B_K}{\beta}}+t,$$
then $\P_{f,g}(\hat T_K \leq t)\leq \beta/2,$ and $\P_{f,g}\left( \hat{\Phi}_{K,\alpha}=0\right)\leq \beta$.
An upper bound for $t$ is finally derived from (\ref{controlq}), which concludes the proof.\\

{\bf Acknowledgments.} We are grateful to the referee, whose discussion and comments allowed us to improve some of our results and to clarify some important points. We also acknowledge the support of the
French Agence Nationale de la Recherche (ANR), under grant ATLAS
(JCJC06-137446) ''From Applications to Theory in Learning and
Adaptive Statistics'' and under grant Calibration
(ANR-2011-BS01-010-02).

\newpage

\section{Supplementary materials}

\subsection{Simulation study}
\subsubsection{Presentation of the simulation study}
In this section, we study our testing procedures from a practical
point of view. We consider $\X=[0,1]$ or $\X=\R$, $n=100$ and $\nu$
the Lebesgue measure on $\X$. $N^1$ and $N^{-1}$ denote two
independent Poisson processes with intensities $f$ and $g$ on $\X$
with respect to $\mu$ with $d\mu=100 d\nu$. We focus on several
couples of
 intensities $(f,g)$ defined on $\X$  and such that $\int_\X f(x) d\nu_x= \int_\X g(x) d\nu_x= 1$. We choose $\alpha =0.05$. \\
Conditionally on the number of points of both processes $N^1$ and
$N^{-1}$, the points of $N^1$ and $N^{-1}$ form two independent
samples of i.i.d. variables  with densities $f$  and $g$ with
respect to $\nu$. Hence, conditionally on the number of points of
$N^1$ and $N^{-1}$, any  test for the classical two-sample
problem for i.i.d. samples can be used here. We compare our tests to the conditional
Kolmogorov-Smirnov test. Thus we consider five testing procedures,
that we respectively denote by KS, Ne, Th, G, E. The testing procedure KS corresponds to the conditional Kolmogorov-Smirnov test. The testing procedures Ne and Th respectively correspond to $\Phi_\alpha^{(1)}$ and $\Phi_\alpha^{(2)}$ defined in \emph{[Multiple kernels case - Example 1]} and
\emph{[Multiple kernels case - Example 2]} with $\bar{J}=7$ and
$\tilde{J}=6$. The testing procedures G and E are similar to the test
$\Phi_\alpha^{(4)}$ defined in \emph{[Multiple kernels case - Example 4]}. For G, we consider the standard
Gaussian approximation kernel defined by
$k(x)=(2\pi)^{-1/2}\exp(-x^2/2)$ for all $x\in \R$ and for E,
 we consider the Epanechnikov approximation kernel defined by $k(x)=(3/4)(1-x^2)\1_{|x|\leq 1}$. For both tests, we take $\{h_m, m\in \M\}=
\{1/24, 1/16, 1/12, 1/8 ,1/4, 1/2\}$ and the corresponding collection
of kernels $\{K_m,m\in \M\}$ given for all $m$ in $\M$ by
$K_m(x,x')= \frac{1}{h_m}k\pa{\frac{x-x'}{h_m}} $. We also take for both
tests $w_m=1/|\M|=1/6$.

Let us recall that our tests reject $(H_0)$ when there exists $m$ in $\M$ such that
$\hat{T}_{K_m}> q^{(N)}_{m,1- u_{\alpha}^{(N)}e^{-w_m}}$
where $N$ is the pooled process obtained from $N^1$ and $N^{-1}$,
and $ u_{\alpha}^{(N)}$ is defined by (\ref{ualpha}).
 Hence, for each observation of the process $N$ whose number of points is denoted by $N_n$, we have to estimate
$  u_{\alpha}^{(N)}$ and the quantiles $q^{(N)}_{m,1-
u_{\alpha}^{(N)}e^{-w_m}}$. These estimations are done by  classical
Monte Carlo methods based on the simulation of $400000$
independent
 samples of size $N_n$ of i.i.d. Rademacher variables (see  Section \ref{MonteCarlo} for the theoretical study of these Monte Carlo methods when single tests are considered). Half of the samples is used to estimate the distribution of each $\hat{T}_{K_m}^\e$. The other half is used to approximate the conditional probabilities
occurring in (\ref{ualpha}). The approximation of $u_{\alpha}^{(N)}$  is
obtained by dichotomy, such that the estimated conditional
probability occurring in (\ref{ualpha}) is less than $\alpha$, but
as close as possible to $\alpha$. By monotony arguments, this is
equivalent to make $u$ varying on a regular grid of $[0,1]$ with
bandwidth $2^{-16}$, and to choose the approximation of $u_{\alpha}^{(N)}$ as the largest
value of the $u$'s on the grid such that the estimated conditional
probabilities in  (\ref{ualpha}) are less than $\alpha$.

\subsubsection{Simulation results}

We first study the probability of first kind error of each test for three common intensities. The first one is the uniform density on $[0,1]$, the second one is the Beta density with parameters $(2,5)$, and the third one is a Laplace density with parameter $7$. Let
\begin{eqnarray*}
 f_1(x)&=&  \1_{[0,1]}(x), \\
 f_{2,2,5}(x)&=&  \frac{x(1-x)^4}{\int_0^1 x(1-x)^4 dx}\1_{[0,1]}(x),\\
 f_{3,7}(x)&=&  \frac{7}{2} e^{-7|x-1/2|}.
\end{eqnarray*}
Taking $f$ as one of these three functions,  we realize $5000$
simulations of two independent Poisson processes $N^1$ and $N^{-1}$
both with intensity $f$ w.r.t. to $\mu$. For each simulation, we
determine the conclusions of the tests KS, Ne, Th, G and E, where
the critical values of our four last tests are approximated by the Monte
Carlo methods described above. The probabilities of first kind error of the tests are estimated
by the number of rejections for these tests divided by $5000$. The
results are given in the following table:

\begin{center}
\begin{tabular}{|c|c|c|c|c|c|}
\hline
$f$ & KS & Ne & Th & G & E \\
\hline
\hline
{\bf ${f}_1$} &$0.053$& $0.049$& $0.045$& $0.053$ & $0.053$ \\
\hline
{\bf ${f}_{2,2,5}$} &$0.053$& $0.047 $& $0.043 $& $0.051 $ & $ 0.050$ \\
\hline
{\bf ${f}_{3,7}$} &$0.0422 $& $ 0.0492 $& $0.0438  $& $0.054 $ & $0.055 $  \\
\hline

\end{tabular}
\end{center}

We then study the probability of second kind error of each test, or more precisely the power of each test, for several alternatives. We consider alternative intensities $(f,g)$ such that $f=f_1$ and $g$ is successively equal to intensities that are classical examples in wavelet settings, and are defined by
\begin{eqnarray*}
 && g_{1,a,\varepsilon}(x)= (1+\varepsilon)\1_{[0,a)}(x) + (1-\varepsilon)\1_{[a,2a)}(x)+\1_{[2a,1)}(x),\\
&& g_{2,\eta}(x)=\pa{1+\eta \sum_j \frac{h_j}{2}(1+\mbox{sgn}(x-p_j))} \frac{\1_{[0,1]}(x)}{C_2(\eta)},\\
&& g_{3,\varepsilon}(x)= (1-\varepsilon) \1_{[0,1 ]}(x) + \varepsilon \pa{\sum_{j} g_j\pa{1+\frac{|x-p_j|}{w_j}}^{-4}}\frac{\1_{[0,1]}(x)}{0.284},\\
\end{eqnarray*}
where $p,h,g,w,\e$ are defined as in \cite{MBP} \footnote{
$$
\begin{tabular}{cccccccccccccc}
p= &(& 0.1&0.13&0.15&0.23&0.25&0.4&0.44& 0.65&0.76& 0.78& 0.81&)\\
h=&(&4&-4&3&-3&5&-5&2&4&-4&2&-3&)\\
g= &(&4&5&3&4&5&4.2&2.1&4.3&3.1&5.1&4.2&) \\
w= &(&0.005&0.005&0.006&0.01&0.01&0.03&0.01&0.01&0.005&0.008&0.005&)
\end{tabular}$$},
$0<\e\leq 1$, $0<a<1/2 $, $\eta >0$ and  $C_2(\eta)$ is such that
$\int_0^1  g_{2,\eta}(x)dx=1$. We also consider alternative intensities $(f,g)$ such that $f$ is equal to the above Laplace density $f_{3,7}$ with parameter $7$, or to the Laplace density $f_{3,10}$ with parameter $10$, such that $f_{3,10}(x)= 5 e^{-10|x-1/2|},$
and $g=g_{4,1/2,1/4}$ is the density of a Gaussian variable with expectation $1/2$ and standard deviation $1/4$.

\smallskip

For each alternative $(f,g)$, we realize $1000$ simulations of two independent Poisson processes $N^1$ and $N^{-1}$  with respective intensities
 $f$ and $g$ w.r.t. $\mu$.   For each simulation, we determine the conclusions of the tests KS, Ne, Th, G and E, where the critical values of our four last tests are still approximated by the Monte Carlo methods described above. The powers of the tests are estimated by the number of
rejections divided by $1000$. The results are
summarized in Figures 1 and 2 where in each column,
 the estimated power is represented as a dot for every test. The triangles represent the upper and lower bounds of an
asymptotic confidence interval with confidence level $99\%$, with
variance estimation.

\begin{figure}[h!]
\label{g1aeg2eta}
\begin{center}
\begin{tabular}{cc}
\includegraphics[scale=0.33]{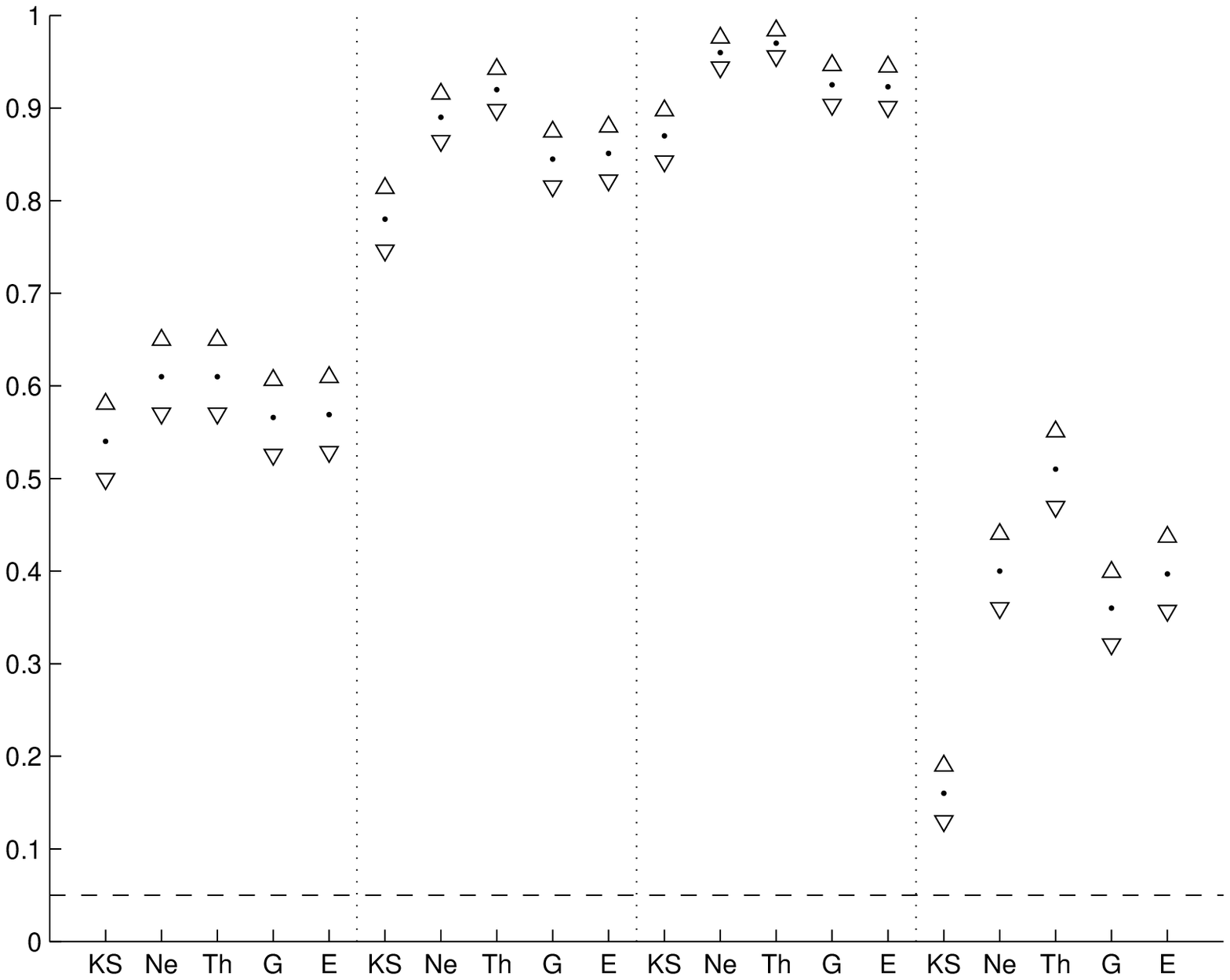} ~~&~~ \includegraphics[scale=0.33]{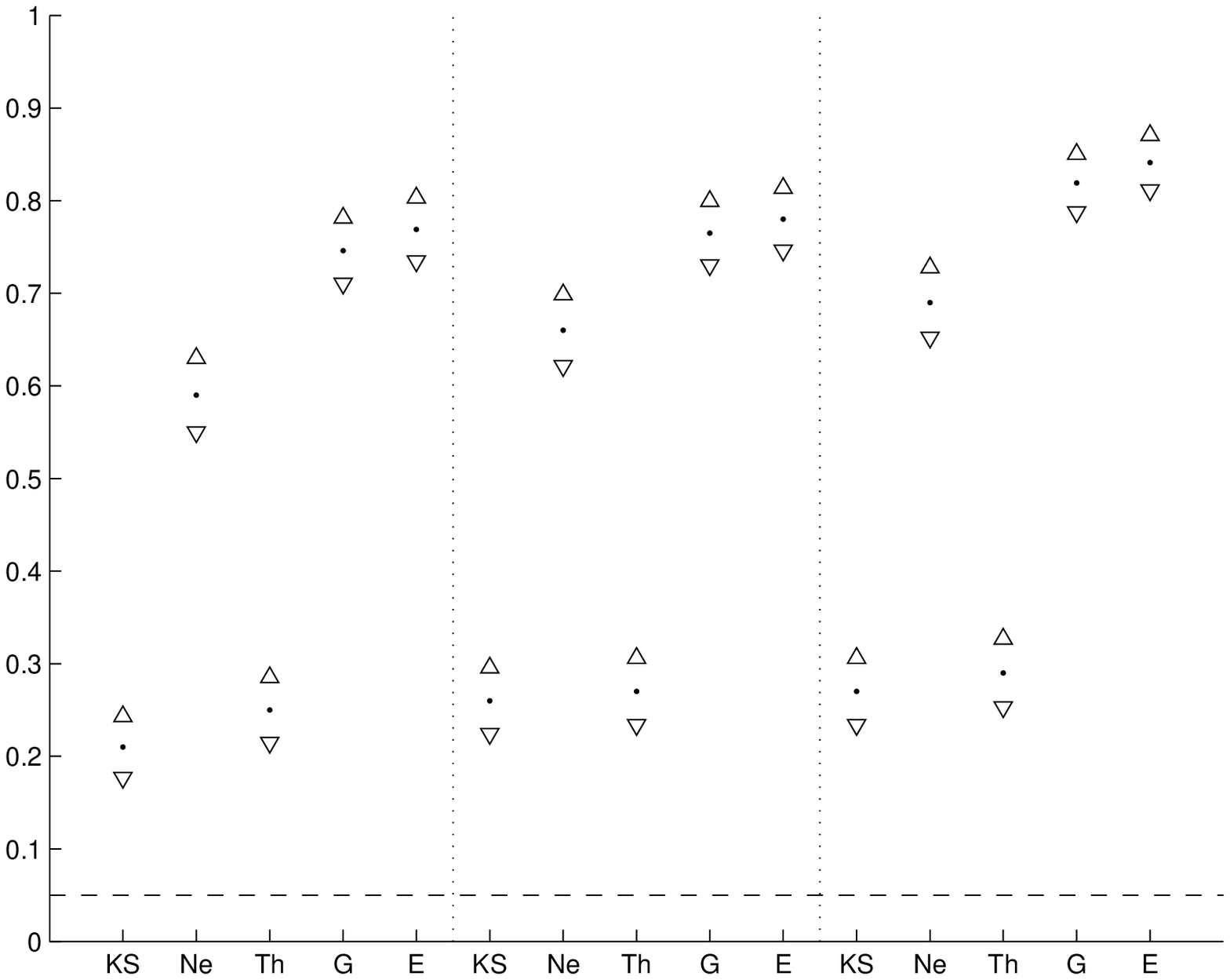}
\end{tabular}
\end{center}
\caption{Left: $(f,g)=(f_1,g_{1,a,\varepsilon})$. Each column corresponds
respectively to $(a,\varepsilon)=(1/4, 0.7),~(1/4,0.9), ~(1/4,1)$ and $(1/8,1)$.
Right: $(f,g)=(f_1,g_{2,\eta})$. Each column corresponds
respectively to $\eta=4,~8$ and $15$.}
\end{figure}

\newpage

\begin{figure}[h!]
\label{g3elaplacegauss}
\begin{center}
\begin{tabular}{cc}
\includegraphics[scale=0.33]{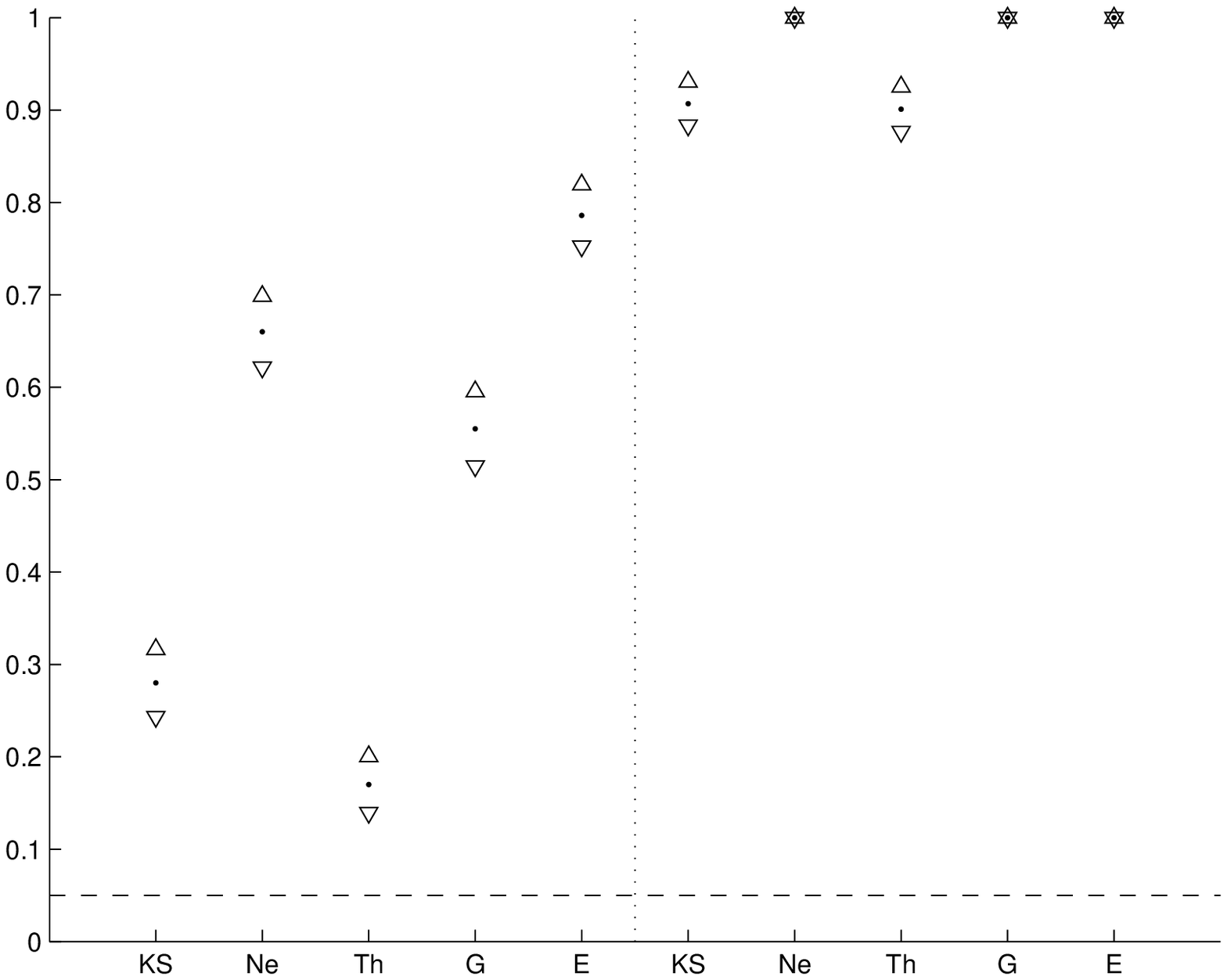}~~&~~\includegraphics[scale=0.33]{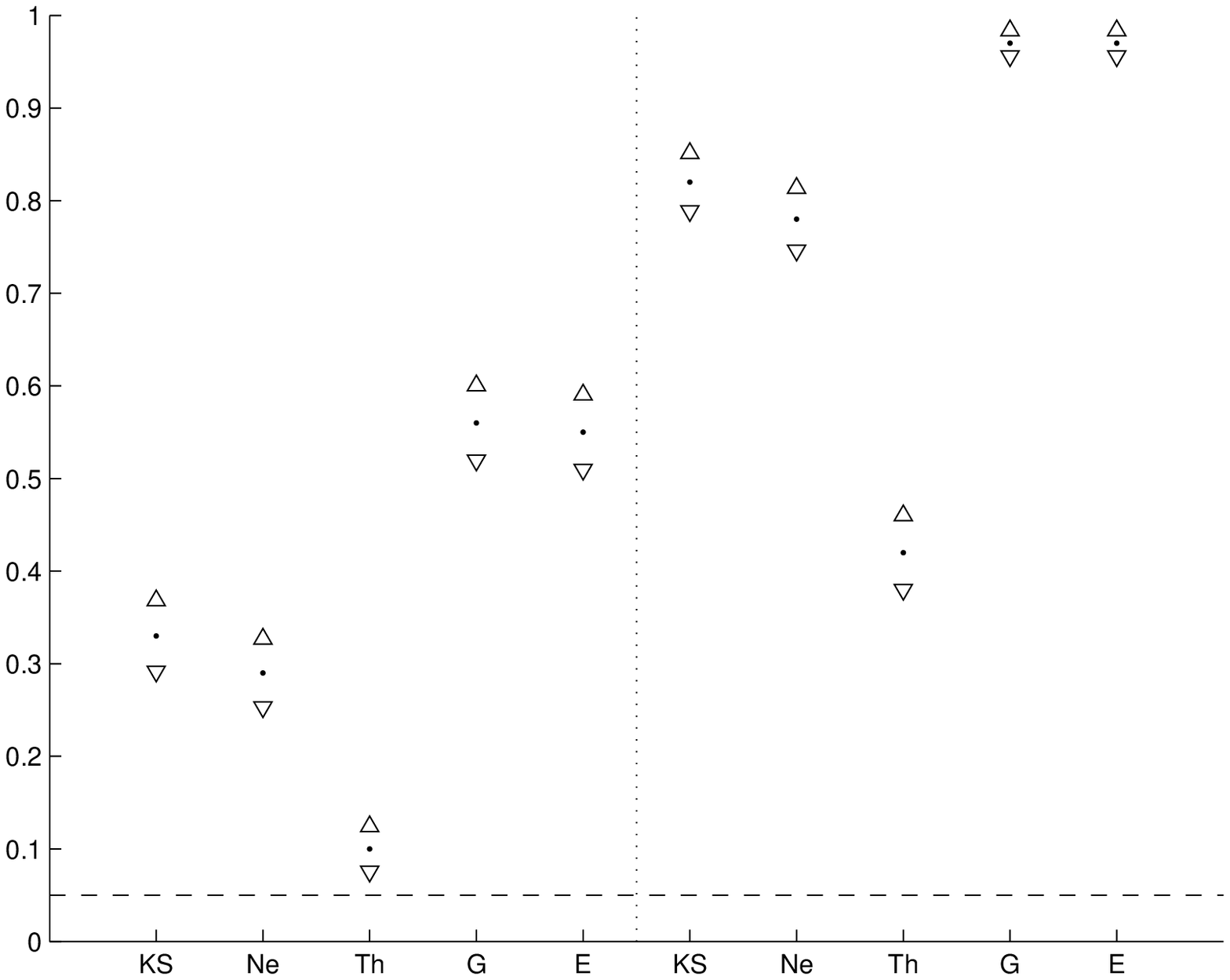}
\end{tabular}
\end{center}
\caption{Left: $(f,g)=(f_1,g_{3,\e})$. The two columns correspond
respectively to $\varepsilon =0.5$ and $1$. Right:
$(f,g)=(f_{3,\lambda}, g_{4,1/2,1/4})$. The two columns correspond
respectively to $\lambda=7$ and $\lambda=10$. }
\end{figure}

In all cases, the tests G and E based on approximation kernels
are more powerful (sometimes even about 4 times more powerful) than
the KS test.  This is also the case for the test Ne, except for the
last example. The test Th is more powerful than the KS test for the
alternatives $(f,g)=(f_1,g_{1,a,\e})$, but it fails to improve the KS
test for the other alternatives. We conjecture that the test Th
consists in the aggregation of too many single tests. We can finally
notice that the test E strongly performs for every considered
alternative, except in a sparse case, where the test  E is less
powerful than the test Th (see Figure 1). Our conclusion is that the
test E is a good practical choice, except maybe when sparse
processes are involved. Aggregating the tests E and Th in such cases
would probably be a good compromise.

\subsection{Proof of Theorem \ref{testmulti} and Theorem \ref{testmultapproxkern}} It is clear from the definition of $u_{\alpha}^{(N)}$ that the test
defined by $\Phi_{\alpha} $ is of level $\alpha$. Obviously, by
Bonferonni's inequality,  $u_{\alpha}^{(N)} \geq \alpha $, hence,
setting $\alpha_m =\alpha e^{-w_m}$, we have
\begin{eqnarray*}
&&\P_{f,g}\pa{ \exists m \in \M, \hat{T}_{K_m} > q^{(N)}_{K_m,1- e^{-w_m} u_{\alpha}^{(N)}} } \\
&&\geq \P_{f,g}\pa{ \exists m \in \M, \hat{T}_{K_m} > q^{(N)}_{K_m,1- \alpha_m} } \\
&&\geq 1- \P_{f,g}\pa{ \forall m \in \M, \hat{T}_{K_m} \leq q^{(N)}_{K_m,1- \alpha_m} } \\
&&\geq  1- \inf_{m \in \M} \P_{f,g}\pa{\hat{T}_{K_m} \leq q^{(N)}_{K_m,1- \alpha_m} } \\
&&\geq  1-\beta,
\end{eqnarray*}
as soon as there exists $m$ in $\M$ such that $
\P_{f,g}\pa{\hat{T}_{K_m} \leq q^{(N)}_{K_m,1- \alpha_m} }\leq
\beta$.
We can now apply Theorem \ref{singlerror}, replacing $\ln(2/\alpha)$ by $(\ln(2/\alpha)+w_m)$, to conclude the proof.

\subsection{Proof of Corollary \ref{vitesses}}
Let us first find an upper bound for\\
 $\rho(\Phi_\alpha^{(1)}, \mathcal{B}_{\delta, \gamma,\infty}(R,R',R'') ,\beta)$.
Considering \eqref{nested}, we in fact only need to find a sharp upper bound for the right hand side of the inequality when $(f,g)$ belongs to $\mathcal{B}_{\delta, \gamma,\infty}(R,R',R'')$.
So let us assume here that $(f,g)\in\mathcal{B}_{\delta, \gamma,\infty}(R,R',R'')$. Then $(f-g)\in \mathcal{B}_{2,\infty}^\delta(R)$, and it is well known (see
\cite{MBP} for instance) that in this case,
$$ \norm{(f-g)-\Pi_{S_J}(f-g)}^2 \leq C(\delta)R^2 2^{-2J\delta}.$$
 Since the constant $C(\alpha,\beta,\norm{f}_\infty,\norm{g}_\infty)$ in \eqref{nested} can be upper bounded by a constant $C(\alpha,\beta,R'')$, the right hand side of \eqref{nested} can be upper bounded by
$$C(\alpha,\beta,\delta,R,R'')\inf_{J\in\M_{\bar{J}}}\left\{2^{-2J\delta}+ (\ln (J+2)) \frac{2^{J/2}}{n}\right\}.$$
Now, taking
$$J^*=\left\lfloor
\log_2\pa{\left(\frac{n}{\ln \ln
n}\right)^{\frac{2}{4\delta+1}}}\right\rfloor,$$
\begin{eqnarray*}
C(\alpha,\beta,\delta,R,R'')&&\inf_{J\in\M_{\bar{J}}}\left\{
 2^{-2J\delta}+ (\ln (J+2)) \frac{2^{J/2}}{n}\right\}\\
&& \leq C(\alpha,\beta,\delta,R,R'')
 \left\{2^{-2J^*\delta}+ (\ln (J^*+2)) \frac{2^{J^*/2}}{n}\right\}\\
 && \leq C(\alpha,\beta,\delta,R,R'')\left(\frac{n}{\ln \ln n}\right)^{-\frac{2\delta}{4\delta+1}}.
 \end{eqnarray*}
 This leads to $$\rho(\Phi_\alpha^{(1)}, \mathcal{B}_{\delta, \gamma,\infty}(R,R',R'') ,\beta)\leq C( \alpha, \beta,\delta,R,R'')
\left(\frac{n}{\ln \ln n}\right)^{-\frac{2\delta}{4\delta+1}}.$$
Of course a similar upper bound applies to $\Psi_\alpha$.

\smallskip

Let us now find an upper bound for $\rho(\Phi_\alpha^{(2)}, \mathcal{B}_{\delta, \gamma,\infty}(R,R',R'') ,\beta)$.
Considering \eqref{thresh}, we only need to find an upper bound for the right hand side of the inequality when $(f,g)$ belongs to $\mathcal{B}_{\delta, \gamma,\infty}(R,R',R'')$.
Let $J$ be an integer that will be chosen later. As in \cite{MBP},
for any $m\subset \Lambda_J=\{(j,k), j\in\{0,\ldots, J-1\},\
k\in\{0,\ldots,2^j-1\}\}$, one can write
$$ \norm{(f-g)-\Pi_{S_m}(f-g)}^2= \norm{(f-g)-\Pi_{S_J}(f-g)}^2+\norm{\Pi_{S_m}(f-g)-\Pi_{S_J}(f-g)}^2.$$
Let us define the coefficients: $\alpha_\lambda=\langle f-g,\p_\lambda\rangle$ for every $\lambda$ in $\{0\}\cup \Lambda_J$, and let us consider $m$ such that $\{\alpha_\lambda, \lambda \in m\}$ is the set
of the $D$ largest coefficients among $\{\alpha_\lambda, \lambda \in
\{0\}\cup \Lambda_J\}$. From \cite{MBP} p.36 for instance, we deduce that
$$\norm{\Pi_{S_m}(f-g)-\Pi_{S_J}(f-g)}^2 \leq C(\gamma)R'^{2+4\gamma} D^{-2\gamma}.$$
As above, we also have that
$$\norm{(f-g)-\Pi_{S_J}(f-g)}^2\leq C(\delta)R^2 2^{-2J\delta}.$$
 Since the constant $C(\alpha,\beta,\norm{f}_\infty,\norm{g}_\infty)$ in \eqref{thresh} can be upper bounded by a constant $C(\alpha,\beta,R'')$, taking
 $$J=\lfloor \log_2
n^\epsilon \rfloor+1$$ for some $\epsilon>0$, the right hand side of
\eqref{thresh} is upper bounded by
$$ C(\alpha,\beta,\delta,\gamma, R,R',R'')\left\{n^{-2\epsilon\delta} +  D^{-2\gamma} + \frac{\epsilon D\ln n}{n}\right\}.$$
Now, taking $D=\lfloor (n/\ln n)^{1/(2\gamma+1)} \rfloor$, and
$\epsilon>\gamma/(\delta(2\gamma+1))$, one obtains that when
$\delta<\gamma/2$, then $D\leq 2^J$, and
$$\rho(\Phi_\alpha^{(2)}, \mathcal{B}_{\delta, \gamma,\infty}(R,R',R'') ,\beta)\leq
 C(\alpha,\beta,\delta,\gamma,R,R',R'') \left(\frac{n}{\ln  n}\right)^{-\frac{\gamma}{1+2\gamma}}.$$
Since this upper bound also applies to  $\Psi_\alpha$, one has
\begin{multline*}
\rho(\Psi_\alpha, \mathcal{B}_{\delta, \gamma,\infty}(R,R',R'')
,\beta)\\
\leq
 C(\alpha,\beta,\delta,\gamma,R,R',R'') \inf\left\{\left(\frac{n}{\ln \ln n}
\right)^{-\frac{2\delta}{4\delta+1}}~,~\left(\frac{n}{\ln
n}\right)^{-\frac{\gamma}{1+2\gamma}}\right\}.
\end{multline*}

\subsection{Proof of the lower bounds}\label{proofLB}
We give here the arguments to derive from the results given in \cite{MBP} lower bounds for
the minimax separation rates over $\mathcal{B}_{\delta, \gamma,\infty}(R,R',R'')$.
As usual, we introduce a finite subset $\mathcal{C}$ of
$\mathcal{B}_{\delta, \gamma,\infty}(R,R',R'')$, composed of couples
of intensities which are particularly difficult to distinguish. Here
one can use the finite subset of possible intensities
$\mathcal{S}_{M,D,r}$ that has been defined in \cite{MBP} Equation
(6.4),
 and define
$$\mathcal{C}=\{(f,g), f=\rho\1_{[0,1]} \mbox{ and } g \in \mathcal{S}_{M,D,r}\},$$
for some fixed positive $\rho$.
Next the computations of the lower bounds of \cite{MBP} can be completely reproduced once we remark that the likelihood ratio
$$\frac{d\P_{\rho{{\bf 1}_{[0,1]}},g}}{d\P_{\rho\1_{[0,1]},\rho\1_{[0,1]}}}(N^1,N^{-1}) =
\frac{d\P_{\rho\1_{[0,1]}}}{d\P_{\rho\1_{[0,1]}}}(N^1)\times\frac{d\P_{g}}{d\P_{\rho\1_{[0,1]}}}(N^{-1}),$$
where on the left hand side $\P_{f,g}$ represents the joint
distribution of two independent Poisson processes $N^1$ and
$N^{-1}$, with respective intensities $f$ and $g$, and on the right
hand side $\P_f$ represents the distribution of one Poisson
process with intensity $f$. This means that the likelihood
ratios that have been considered in \cite{MBP} are exactly the ones
we need here to compute the expected lower bounds. The results are
consequently identical.

\subsection{Proof of Corollary \ref{vitessessobolev}}

Considering \eqref{oraclenoyausob}, we mainly have to find a sharp upper bound for \begin{multline*}
\inf_{(m_1,m_2)\in \M}
\Bigg\{\norm{(f-g)-k_{m_1,h_{m_2}}*(f-g)}^2\\
+\frac{w_{(m_1,m_2)}}{n}\sqrt{\frac{\norm{f+g}_\infty\norm{f+g}_1\norm{k_{m_1}}^2}{
2^{-dm_2}}}\Bigg\},
\end{multline*}
when $(f,g)$ belongs to $\mathcal{S}_{d}^{\delta}(R,R',R'')$.\\
Let
us first control the bias term
$\norm{(f-g)-k_{m_1,h_{m_2}}*(f-g)}^2$. Plancherel's theorem gives
that when $(f-g)\in\mathbb{L}^1(\R^d)\cap\mathbb{L}^2(\R^d)$,
\begin{eqnarray*}
(2\pi)^d\norm{(f-g)&-&k_{m_1,h_{m_2}}*(f-g)}^2\\
&&=\norm{(1-\widehat{k_{m_1,h_{m_2}}})\widehat{(f-g)}}^2\\
&&=\int_{\R^d} \left|1-\widehat{k_{m_1}}(2^{-m_2}u)\right|^2(u)
(\widehat{f-g})^2(u) d\nu_u.
\end{eqnarray*}
Assume now that $(f,g)\in\mathcal{S}_{d}^{\delta}(R,R',R''),$ and
take $m_1^*=\textrm{min}\{m_1\in\M_1, m_1 \geq \delta\}$. Note that
since $\norm{\widehat{k_{m_1^*}}}_\infty<+\infty$ and
$\widehat{k_{m_1^*}}$ satisfies the condition (\ref{nonintkernel}),
there also exists some constant $C(\delta)>0$ such that
$$\mbox{Ess sup}_{u\in\R^d\setminus \{0\}}\frac{|1-
\widehat{k_{m_1^*}}(u)|}{\norm{u}_d^{\delta}}\leq C(\delta).$$ Then
$$\norm{(f-g)-k_{m_1^*,h_{m_2}}*(f-g)}^2\leq \frac{C(\delta)}{(2\pi)^d} \int_{\R^d}
\norm{2^{-m_2}u}_d^{2\delta}(\widehat{f-g})^2(u) d\nu_u,$$ and
since $(f-g)\in\mathcal{S}_{d}^{\delta}(R)$,
$$\norm{(f-g)-k_{m_1^*,h_{m_2}}*(f-g)}^2\leq 2^{-2\delta
m_2}C(\delta) R^2.$$
Furthermore, $\norm{k_{m_1^*}}^2\leq C(\delta),$
so
\begin{multline*}
\inf_{(m_1,m_2)\in \M}
\Bigg\{\norm{(f-g)-k_{m_1,h_{m_2}}*(f-g)}^2\\
+\frac{w_{(m_1,m_2)}}{n}\sqrt{\frac{\norm{f+g}_\infty\norm{f+g}_1\norm{k_{m_1}}^2}{
2^{-dm_2}}}\Bigg\}\\
\leq C(\delta,\alpha,\beta,R)\inf_{m_2\in \M_2} \Bigg\{2^{-2\delta
m_2}+\frac{w_{(m_1^*,m_2)}}{n}\sqrt{\frac{\norm{f+g}_\infty\norm{f+g}_1}{2^{-dm_2}}}\Bigg\}.
\end{multline*}
Choosing $$m_2^*=\left\lfloor \log_2 \left( \left(\frac{n}{\ln \ln
n}\right)^{\frac{2}{d+4\delta}}\right)\right\rfloor$$ leads to
$$2^{-2\delta m_2^*}\leq 2^{2\delta}\left(\frac{\ln \ln
n}{n}\right)^{\frac{4\delta}{d+4\delta}},$$ and since
$w_{(m_1^*,m_2^*)}\leq C(\delta,d)\ln \ln n$,
$$\frac{w_{(m_1^*,m_2^*)}}{n}\sqrt{\frac{\norm{f+g}_\infty\norm{f+g}_1}{
2^{-d m_2^*}}} \leq
C(\delta,d)\sqrt{\norm{f+g}_\infty\norm{f+g}_1}\left(\frac{\ln\ln
n}{n}\right)^{\frac{4\delta}{d+4\delta}}.$$ Noting that
$1/n\leq \left(\ln\ln
n/n\right)^{4\delta/(d+4\delta)},$ when
$(f,g)\in\mathcal{S}_d^\delta(R,R',R'')$,
\begin{multline*}
C(\alpha,\beta)\Bigg(\inf_{(m_1,m_2)\in \M}
\Bigg\{\norm{(f-g)-k_{m_1,h_{m_2}}*(f-g)}^2\\
+\frac{w_{(m_1,m_2)}}{n}\sqrt{\frac{\norm{f+g}_\infty\norm{f+g}_1\norm{k_{m_1}}^2}{
2^{-dm_2}}}\Bigg\}+\frac{\norm{f+g}_\infty}{n }\Bigg)\\
\leq C(\delta,\alpha,\beta,R,R',R'',d)\left(\frac{\ln\ln
n}{n}\right)^{\frac{4\delta}{d+4\delta}}.
\end{multline*}
This concludes the proof of Corollary \ref{vitessessobolev}.

\subsection{Proof of Corollary \ref{vitessesanisotrop}}

As in the previous section, considering
\eqref{oraclenoyaunik}, we here have to find a sharp upper bound
for
$$ \inf_{m_2\in \M_2}
\Bigg\{\norm{(f-g)-k_{1,h_{m_2}}*(f-g)}^2
+\frac{w_{(1,m_2)}}{n}\sqrt{\frac{\norm{f+g}_\infty\norm{f+g}_1\norm{k_1}^2}{
\prod_{i=1}^dh_{m_2,i}}}\Bigg\}.$$ Let us first evaluate
$\norm{(f-g)-k_{1,h}*(f-g)}^2$ when
$(f-g)\in\mathcal{N}_{2,d}^{\delta}(R)$, and $h=(h_1,\ldots,h_d)$.\\
For $x=(x_1,\ldots,x_d)\in\R^d$, let
$b(x)=k_{1,h}*(f-g)(x)-(f-g)(x)$. Then
$$b(x)=\int_{\R^d}
k_1(u_1,\ldots,u_d)(f-g)(x_1+u_1h_1,\ldots,x_d+u_dh_d)du_1\ldots
du_d-(f-g)(x),$$ and since $\int_{\R^d}
k_1(u_1,\ldots,u_d)du_1\ldots du_d=1$,

\begin{eqnarray*}
b(x)&=&\int_{\R^d}
k_1(u_1,\ldots,u_d)\Big[(f-g)(x_1+u_1h_1,\ldots,x_d+u_dh_d)\\
&&-(f-g)(x_1,\ldots,x_d)\Big]du_1\ldots
du_d\\
&= &\sum_{i=1}^d b_i(x),
\end{eqnarray*}
where for $i=1\ldots d$,

\begin{multline*} b_i(x)=\int_{\R^d} k_1(u_1,\ldots,u_d)\Big[(f-g)(x_1+u_1h_1,\ldots,x_i+u_ih_i,x_{i+1},\ldots,x_d)\\
-(f-g)(x_1+u_1h_1,\ldots,x_i,x_{i+1},\ldots,x_d)\Big]du_1\ldots
du_d,\end{multline*}
 As in the proof of Proposition 1.5 p. 13 of
\cite{Sacha}, using the Taylor expansion of $(f-g)$ in the $i$th
direction and the fact that $\int_\R k_{1,i}(u_i)u_i^jdu_i=0$ for
$j=1\ldots \Delta_i$, we obtain that

\begin{eqnarray*}
b_i(x)&=&\int_{\R^d} k_1(u)\frac{(u_ih_i)^{\lfloor
\delta_i\rfloor}}{(\lfloor \delta_i\rfloor\!-\!1)!}
\Bigg[\int_0^1(1-\tau)^{\lfloor \delta_i\rfloor-1}\\
&&D_i^{\lfloor \delta_i\rfloor}(f-g)(x_1+u_1h_1,\ldots,x_i+\tau
u_ih_i,x_{i+1},\ldots,x_d)d\tau\Bigg]
du,
\end{eqnarray*}
So,
\begin{eqnarray*}
b_i(x)&=&\int_{\R^d} k_1(u)\frac{(u_ih_i)^{\lfloor
\delta_i\rfloor}}{(\lfloor
\delta_i\rfloor\!-\!1)!}\Bigg[\int_0^1(1-\tau)^{\lfloor
\delta_i\rfloor-1}\\
&&\Big(D_i^{\lfloor \delta_i\rfloor}(f-g)(x_1+u_1h_1,\ldots,x_i+\tau
u_ih_i,x_{i+1},\ldots,x_d)\\
&&-D_i^{\lfloor
\delta_i\rfloor}(f-g)(x_1+u_1h_1,\ldots,x_i,\ldots,x_d)\Big)d\tau\Bigg]du.
\end{eqnarray*}
Hence, by using twice Lemma 1.1 p. 13 of \cite{Sacha} extended to
the spaces $\R^d\times\R$ and $\R^d\times\R^d$,

\begin{eqnarray*}
\norm{b_i}_2^2&\leq & \int_{\R^d} \Bigg(\int_{\R^d}
|k_1(u)|\frac{|u_ih_i|^{\lfloor \delta_i\rfloor}}{(\lfloor
\delta_i\rfloor\!-\!1)!}\Bigg[\int_0^1(1-\tau)^{\lfloor
\delta_i\rfloor-1}\\
&&\Big|D_i^{\lfloor \delta_i\rfloor}(f-g)(x_1+u_1h_1,\ldots,x_i+\tau
u_ih_i,x_{i+1},\ldots,x_d)\\
&&-D_i^{\lfloor
\delta_i\rfloor}(f-g)(x_1+u_1h_1,\ldots,x_i\ldots,x_d)\Big|d\tau\Bigg]du\Bigg)^2dx\\
&\leq & \Bigg[\int_{\R^d} |k_1(u)|\frac{|u_ih_i|^{\lfloor
\delta_i\rfloor}}{(\lfloor \delta_i\rfloor\!-\!1)!} \Bigg(
\int_{\R^d}\Big[\int_0^1(1-\tau)^{\lfloor
\delta_i\rfloor-1}\\
&&\Big|D_i^{\lfloor \delta_i\rfloor}(f-g)(x_1+u_1h_1,\ldots,x_i+\tau
u_ih_i,x_{i+1},\ldots,x_d)\\
&&-D_i^{\lfloor
\delta_i\rfloor}(f-g)(x_1+u_1h_1,\ldots,x_i\ldots,x_d)\Big|d\tau\Big]^2dx\Bigg)^{1/2}du\Bigg]^2
\end{eqnarray*}
 and
\begin{eqnarray*}
\norm{b_i}_2^2&\leq & \Bigg[\int_{\R^d}
|k_1(u)|\frac{|u_ih_i|^{\lfloor \delta_i\rfloor}}{(\lfloor
\delta_i\rfloor\!-\!1)!} \Bigg( \int_0^1(1-\tau)^{\lfloor
\delta_i\rfloor-1}\Big(\int_{\R^d}\\
&&\Big|D_i^{\lfloor \delta_i\rfloor}(f-g)(x_1+u_1h_1,\ldots,x_i+\tau
u_ih_i,x_{i+1},\ldots,x_d)\\
&&-D_i^{\lfloor
\delta_i\rfloor}(f-g)(x_1+u_1h_1,\ldots,x_i\ldots,x_d)\Big|^2dx\Big)^{1/2}d\tau\Bigg)
du\Bigg]^2.
\end{eqnarray*}
When $(f-g)\in\mathcal{N}_{2,d}^{\delta}(R)$,
$$\norm{b_i}_2\leq C(\delta_i)R \int_{\R^d} |k_1(u)| |u_ih_i|^{\delta_i}
du\leq C(\delta_i)\left(\int_{\R^d} |k_1(u)| |u_i|^{\delta_i}
du\right) Rh_i^{\delta_i}.$$
So,
$$\norm{k_{1,h}*(f-g)-(f-g)}\leq
C(\delta)R\sum_{i=1}^dh_i^{\delta_i}.$$
Let us now find some $m_2$ in $\M_2$ giving a sharp upper bound for
$$
\inf_{m_2\in \M_2} \Bigg\{\norm{(f-g)-k_{1,h_{m_2}}*(f-g)}^2
+\frac{w_{(1,m_2)}}{n}\sqrt{\frac{\norm{f+g}_\infty\norm{f+g}_1\norm{k_1}^2}{
\prod_{i=1}^dh_{m_2,i}}}\Bigg\}.
$$
Let ${1}/{\bar{\delta}}=\sum_{i=1}^d {1}/{\delta_i},$ and
choose $m_2^*=(m_{2,1}^*,\ldots,m_{2,d}^*)$ in $\M_2$, with
$$m_{2,i}^*=\left\lfloor \log_2 \left(
\left(\frac{n}{\ln \ln
n}\right)^{\frac{2\bar{\delta}}{\delta_i(1+4\bar{\delta})}}\right)\right\rfloor,$$
for every $i=1\ldots d$. Since
$h_{m_2^*}=(2^{-m_{2,1}^*},\ldots,2^{-m_{2,d}^*})$,
$$\norm{(f-g)-k_{1,h_{m_2^*}}*(f-g)} \leq
C(\delta,R)\sum_{i=1}^d2^{-m_{2,i}^*\delta_i},$$ so
$$\norm{(f-g)-k_{1,h_{m_2^*}}*(f-g)}^2\leq C(\delta,R)d^2
\left(\frac{\ln\ln
n}{n}\right)^{\frac{4\bar{\delta}}{1+4\bar{\delta}}}.$$
 Moreover, it is easy to see that
$w_{(1,m_2^*)}\leq C(\delta,d)\ln \ln n,$ and hence
\begin{eqnarray*}
&&\frac{w_{(1,m_2^*)}}{n}\sqrt{\frac{\norm{f+g}_\infty\norm{f+g}_1\norm{k_1}^2}{
\prod_{i=1}^d 2^{-m_{2,i}^*}}} \\
&&\leq C(\delta,\alpha,\beta,R',R'',d)\frac{\ln\ln
n}{n}\left(\frac{n}{\ln\ln
n}\right)^{\sum_{i=1}^d\frac{\bar{\delta}}{(1+4\bar{\delta})\delta_i}}\\
&&\leq C(\delta,\alpha,\beta,R',R'',d)\left(\frac{\ln\ln
n}{n}\right)^{\frac{4\bar{\delta}}{1+4\bar{\delta}}}.
\end{eqnarray*}
Since $$\frac{1}{n}\leq \left(\frac{\ln\ln
n}{n}\right)^{\frac{4\bar{\delta}}{1+4\bar{\delta}}},$$ when $\ln \ln n\geq 1$,
this ends
the proof of Corollary \ref{vitessesanisotrop}.\\

\end{document}